\documentclass[12p]{article}
\usepackage{amsmath,amssymb,latexsym,theorem,bbm,epsfig,natbib,easybmat,
  eufrak,color}

\newcommand{\Var}{\mathrm{Var}}

\newcommand{\tr}{\mathrm{tr}}
\newcommand{\mspe}{\mathrm{MSPE}}
\newcommand{\imspe}{\mathrm{IMSPE}}
\newcommand{\ent}{\mathrm{Ent}}

\numberwithin{equation}{section}

\newcommand{\proofend}{\hfill$\square$}

\newtheorem{thm}{Theorem}[section]

\theorembodyfont{\rm}
\newtheorem{rem}[thm]{Remark}
\newtheorem{ex}[thm]{Example}

\begin{document}

\title{Optimal designs for the methane flux in troposphere}

\author{S\'andor Baran, Kinga Sikolya\\
Faculty of Informatics,
University of Debrecen, Hungary,\\
\and
Milan Stehl\'\i k
\\
Institut f\"ur Angewandte
Statistik, Johannes Kepler University in Linz, Austria\\
Departamento de Matem\'atica,
Universidad T\'ecnica Federico Santa Mar\'{\i}a,
Valpara\'iso, Chile}

\date{}
\maketitle

\begin{abstract}

The understanding of methane emission and methane absorption plays a
central role both in the atmosphere and on the surface of the Earth.
Several important ecological processes, e.g. ebullition of methane and
its natural microergodicity request better designs for observations in
order to decrease variability in parameter estimation.
Thus, a crucial fact, before the measurements are taken, is to give an
optimal design of the sites where observations should be collected in
order to stabilize the
variability of estimators. In this paper we introduce a realistic
parametric model of covariance and provide theoretical and numerical
results on optimal designs. For parameter estimation D-optimality,
while for prediction integrated mean square error and entropy criteria
are used. We illustrate applicability of obtained benchmark designs
for increasing/measuring the efficiency of the engineering designs for
estimation of methane rate in various temperature ranges and under
different correlation parameters. We show that in most situations
these benchmark designs have higher efficiency.

\par

\par

\end{abstract}

\fontsize{10.95}{14pt plus.8pt minus .6pt}\selectfont

\section{Introduction}
  \label{sec:sec1}
The understanding of methane emission and methane absorption plays a
central role both in the atmosphere (for troposphere see,
e.g., \citet{Vag91}) and on the surface of the Earth (see,
e.g., \citet{Li} regarding the
 methane emissions from natural wetlands and references therein or
 \citet{Jordanova} for efficient and robust model of the methane
 emission from sedge-grass marsh in South Bohemia).
Several important ecological processes, e.g. ebullition of methane and
its natural microergodicity request better designs for observations in
order to decrease variability in parameter estimation
\citep{Jordanova2}. In this context by a
  design we mean a set of locations where
  the investigated process is observed.
Thus, a crucial fact, before the measurements are taken, is to give an
optimal design of the sites where observations should be collected.
\citet{Chemo12} provided a comparison of filling and D-optimal designs
for a one-dimensional design variable, e.g., temperature.
However, such a model oversimplifies the important fact that
variation of other variables, e.g., rates $k_1$ of the considered
modified Arrhenius model,
could disturb the efficiency of the learning process. The latter
statement is
also in agreement with common sense in physical chemistry. In this
paper the difficulties of modelling and design are treated,
mainly by  allowing an Ornstein-Uhlenbeck (OU) sheet error model.

We concentrate on  efficient estimation of the parameters of the
modified Arrhenius model (model  popular  in chemical kinetics), which is
used  by \citet{Vag91} as a flux model of methane in troposphere.
This generalized exponential (GE) model can be expressed as
\begin{equation} \label{eq:GEmodel}
Y=A x^{\mu} {\mathrm e}^{-B x}+\epsilon=\eta(x,\mu,B)+\epsilon,
\end{equation}
where $A, B, \mu\in{\mathbb R}$, $A, \ B\geq 0$,  are constants
and $\varepsilon$ is a random error term.  In the case of correlated  errors
such a model was studied by \citet{Chemo12}, however, in that work
error structures were univariate stochastic processes. For the
case of uncorrelated errors see Bayesian approach of \citet{DeSp94}
and also the work of
\citet{RodSan09} for different optimality criteria and restrictions on
the design space.
In \citet{RodSan09} and  \citet{Chemo12} the authors concentrated on the
Modified-Arrhenius (MA) model, which is equivalent to the GE model
through the change of variable $x=1/t$. This model is useful for
chemical kinetic (mainly
because it is a generalization of Arrhenius model describing the
influence of   temperature $t$ on the rates of chemical
processes, see,
e.g., \citet{Lai84} for general discussion and \citet{LicesioJesus} for
optimal designs).
However, for specific setups, for instance, long temperature ranges,
Arrhenius model is insufficient and
 the Modified Arrhenius (or GE model) appears to be  the
good alternative \citep[see, for instance,][]{Gie97}.
Other applications of model \eqref{eq:GEmodel} in chemistry are related to the
transition state theory (TST) of chemical reactions \citep{TST}.

In practical chemical kinetics  two steps are taken:
first the rates $k_1$ are estimated (typically with symmetric
estimated error) and then modified Arrhenius model is fitted to the
rates, i.e.,
\begin{equation}
  \label{k1fitting}
k_1= A (1/t)^{\mu} {\mathrm e}^{-B/t}+\widetilde\varepsilon (t).
\end{equation}
Statistically
correct would be to assess both steps by one optimal
experimental planning.
\citet{Chemo12} concentrated on the
second phase, i.e. what is the optimal distribution of
temperature for obtaining
statistically efficient estimators of trend parameters $A, B, \ \mu$
and correlation parameters of the error term $\widetilde\varepsilon$.
In this paper we provide designs both for rates and temperatures,
and in this way substantially generalize the previously studied model.

Correlation is the natural dependence measure
fitting for elliptically symmetric distributions (e.g., Gaussian).
By taking  $s$ (this variable can play, for example, the role of
    atmospheric pressure,
latitude or location of the measuring balloon in
troposphere, either vertically or horizontally) and  temperature $t$ to be
variables of covariance, our model \eqref{eq:GEmodel}  takes a form  of a
stationary process
\begin{equation}
   \label{model}
Y(s,t) = k_1 +\varepsilon (s,t),
\end{equation}
where the design points are taken from a compact design space
$\mathcal{X}=[a_1,b_1]\times [a_2, b_2]$, with $b_1>a_1$ and
$b_2>a_2$, and $\varepsilon (s,t), \ s,t\in {\mathbb R}$, is a
stationary OU
sheet, that is a zero mean Gaussian process with covariance structure
\begin{equation}
   \label{oucov}
{\mathsf E}\,\varepsilon(s_1,t_1)\varepsilon(s_2,t_2)=
\frac{{\widetilde\sigma}^2}{4\alpha\beta}\exp\big
(-\alpha|s_1-s_2|-\beta|t_1-t_2|\big ),
\end{equation}
where $\alpha>0, \ \beta>0, \ \widetilde\sigma>0$. We remark that
$\varepsilon(s,t)$ can also be represented as
\begin{equation*}
\varepsilon (s,t)=\frac{\widetilde\sigma}{2\sqrt{\alpha\beta}}{\mathrm
  e}^{-\alpha s-\beta t}{\mathcal W}\big({\mathrm e}^{2\alpha s},
  {\mathrm e}^{2\beta t}\big),
\end{equation*}
where ${\mathcal W}(s,t), \ s,t\in {\mathbb R}$, is a  standard
  Brownian sheet \citep{Baran, bs}. Covariance structure \eqref{oucov}
  implies that for ${\mathbf d}=(d,\delta), \ d\geq 0, \ \delta \geq 0$,
  the variogram $2\gamma({\mathbf d}):=\Var \big (\varepsilon
  (s+d,t$ $+\delta)-\varepsilon (s,t)\big )$ equals
\begin{equation*}
2\gamma({\mathbf d})=\frac{{\widetilde\sigma}^2}{2\alpha\beta}\Big(1- {\mathrm
  e}^{-\alpha d-\beta \delta}\Big)
\end{equation*}
and the correlation between two measurements depends
on the distance through the semivariogram  $\gamma({\mathbf d})$.

As can be visible from relation \eqref{k1fitting} between rates and
parameters $A,\mu$ and $ B$ of the modified Arrhenius model, the
second variable $s$  
is missing from trend since it is not chemically understood as driving
mechanism of chemical kinetics, however, in this context it is an
environment variable. 

In order to apply the usual notations of spatial modeling \citep{KS} we
introduce $\sigma:=\widetilde\sigma/(2\sqrt{\alpha\beta})$ and instead
of \eqref{oucov} we investigate
\begin{equation}
   \label{oucovmod}
{\mathsf E}\,\varepsilon(s_1,t_1)\varepsilon(s_2,t_2)=
\sigma^2\exp\big (-\alpha|s_1-s_2|-\beta|t_1-t_2|\big ),
\end{equation}
where $\sigma $ is considered as a nuisance parameter.
For discussion on the identifiability of the covariance
parameters see, e.g., \citet{MullerStehlik2}. 

\section{Benchmarking grid designs for the OU sheet with 
constant  trend}
   \label{sec:sec2}

In this section we derive several optimal design results for the case
of constant trend  and regular grids resulting in a Kronecker
product covariance structure. These theoretical contributions
will serve as benchmarks for optimal designs in a methane flux model.
Thus we consider the stationary process
\begin{equation}
   \label{OUmodel}
Y(s,t) = \theta +\varepsilon (s,t)
\end{equation}
with the design points  taken from a compact design space
$\mathcal{X}=[a_1,b_1]\times [a_2, b_2]$, where $b_1>a_1$ and
$b_2>a_2$ and $\varepsilon (s,t), \ s,t\in {\mathbb R}$, is a
stationary Ornstein-Uhlenbeck
sheet, i.e., a zero mean Gaussian process with covariance structure
(\ref{oucovmod}).

\subsection{D-optimality}
   \label{subsec:subsec2.1}
As a first step we derive D-optimal designs, that is arrangements of
design points that maximize the objective function $\Phi(M):=\det
(M)$, where $M$ is the Fisher information matrix of observations of
the random field $Y$. This method, "plugged" from the widely
developed uncorrelated setup, is offering considerable potential for
automatic implementation, although further development is needed
before it can be applied routinely in practice. Theoretical
justifications of using the Fisher information for D-optimal
designing under correlation can be found in \citet{Abt, Pazman07} and
\citet{Stehlik07}.

We investigate grid designs of the form $\big\{(s_i,t_j): \ i=1,2, \ldots
,n, \ j=1,2, \ldots ,m\big\}\subset \mathcal{X}=[a_1,b_1]\times [a_2,
b_2]$, $n,m\geq 2$, and without loss of generality we may assume 
$a_1\leq s_1<s_2< \ldots 
< s_n\leq b_1$ and $a_2\leq t_1<t_2< \ldots < t_m\leq b_2$. Usually,
the grid containing the design points can be arranged arbitrary in the
design space $\mathcal X$, but we also consider restricted
D-optimality, when  $s_1=a_1, \ s_n=b_1$ and $t_1=a_2, \ t_m=b_2$,
i.e. the vertices of ${\mathcal X}$ are included in all  designs.

\subsubsection{Estimation of trend parameter only}
Let us assume first that parameters $\alpha, \beta$ and $\sigma$ of
the covariance structure \eqref{oucovmod} of the OU sheet
$\varepsilon$ are given and we are interested in estimation of the
trend parameter $\theta$. In this case the Fisher information on
$\theta$ based on observations $\big\{Y(s_i,t_j), \ i=1,2,\ldots,
n, \ j=1,2,\ldots, m\big\}$ equals $M_{\theta}(n,m)={\mathbf
  1}^{\top}_{nm}C^{-1}(n,m,r){\mathbf 1}_{nm}$, where ${\mathbf
  1}_k, \ k\in{\mathbb N},$ denotes the column vector of ones of length
$k$,  $r=(\alpha,\beta)^{\top}$, and $C(n,m,r)$ is the covariance matrix
of the observations  \citep{Pazman07, Xia}.  Further, let
$d_i:=s_{i+1}-s_i, \ \ i=1,2, \ldots ,n-1$, and
$\delta_j:=t_{j+1}-t_j, \ j=1,2, \ldots, m-1$,
be the directional distances between two adjacent design points.
With the help of this representation one can prove the following theorem.

\begin{thm}
   \label{trend} \
Consider the OU model \eqref{OUmodel} with covariance structure
\eqref{oucovmod} observed in points $\big\{ (s_i,t_j),
  \ i=1,2,\ldots ,n, \ j=1,2,\ldots ,m\big \}$ and assume that the only
  parameter of interest is the trend parameter $\theta$. In this case
\begin{equation}
  \label{eq:eq3.1}
M_{\theta}(n,m)=\Bigg(1+\sum_{i=1}^{n-1}\frac{1-p_i}{1+p_i}\Bigg)\Bigg(
1+\sum_{j=1}^{m-1}\frac{1-q_j}{1+q_j}\Bigg),
\end{equation}
where $p_i:=\exp(-\alpha d_i)$, $q_j:=\exp(-\beta
\delta_j), \ i=1,2, \ldots ,n-1, \ j=1,2, \ldots ,m-1$, and
  the directionally equidistant design $d_1=d_2=\ldots =d_{n-1}$ and
$\delta_1=\delta _2=\ldots =\delta_{m-1}$  is optimal
for estimation of  $\theta$.
\end{thm}

\subsubsection{Estimation of covariance parameters only}

Assume now that we are interested only in the estimation of the
parameters $\alpha $ and $\beta$ of the Ornstein-Uhlenbeck
sheet. According to the results of \citet{Pazman07} and \citet{Xia}
the Fisher information matrix on $r=(\alpha,\beta)^{\top}$ has the
form
\begin{equation}
  \label{eq:eq4.1}
M_r(n,m)=\begin{bmatrix}
        M_{\alpha}(n,m) &  M_{\alpha,\beta}(n,m) \\
         M_{\alpha, \beta}(n,m) &  M_{\beta }(n,m)
       \end{bmatrix},
\end{equation}
where
\begin{align*}
M_{\alpha}(n,m)&:=\frac 12 \tr \left\{C^{-1}(n,m,r)\frac{\partial
    C(n,m,r)}{\partial
    \alpha}C^{-1}(n,m,r)\frac{\partial C(n,m,r)}{\partial \alpha} \right\}, \\
M_{\beta}(n,m)&:=\frac 12 \tr \left\{C^{-1}(n,m,r)\frac{\partial
    C(n,m,r)}{\partial
    \beta }C^{-1}(n,m,r)\frac{\partial C(n,m,r)}{\partial \beta} \right\}, \\
M_{\alpha,\beta}(n,m)&:=\frac 12 \tr \left\{C^{-1}(n,m,r)\frac{\partial
    C(n,m,r)}{\partial
    \alpha}C^{-1}(n,m,r)\frac{\partial C(n,m,r)}{\partial \beta} \right\},
\end{align*}
and  $C(n,m,r)$ is the covariance matrix of the observations $\big
\{Y(s_i,t_j), \ i\!=\! 1,2,\ldots, n, \ j\!=\! 1,2,\ldots ,m
\big\}$. Note that here
$ M_{\alpha}(n,m)$ and $ M_{\beta}(n,m)$ are Fisher information on
parameters $\alpha$ and $\beta$, respectively, taking the other
parameter as a nuisance.

The following theorem gives the exact form of \ $M_r(n,m)$ \ for the
model \eqref{OUmodel}.

\begin{thm}
    \label{Mrn} \
Consider the OU model \eqref{OUmodel} with covariance structure
\eqref{oucovmod} observed in points $\big\{ (s_i,t_j),
  \ i=1,2,\ldots ,n, \ j=1,2,\ldots ,m\big \}$. Then
\begin{align}
   \label{eq:eq4.2}
M_{\alpha}(n,m)=&\,m\sum_{i=1}^{n-1}\frac{d_i^2p_i^2(1+p_i^2)}{(1-p_i^2)^2},
 \qquad
M_{\beta}(n,m)=n\sum_{j=1}^{m-1}\frac{\delta_j^2q_j^2(1+q_j^2)}{(1-q_j^2)^2},\\
&M_{\alpha,\beta}(n,m)=2\bigg(\sum_{i=1}^{n-1}\frac{d_ip_i^2}{1-p_i^2}\bigg)
\bigg(\sum_{j=1}^{m-1}\frac{\delta_jq_j^2}{1-q_j^2}\bigg), \nonumber
\end{align}
where $d_i, \delta_j$ and $p_i,q_j$ denote the same quantities as
before, i.e. $d_i:=s_{i+1}-s_i, \ \delta_j:=t_{j+1}-t_j$ and
$p_i:=\exp(-\alpha d_i)$, $q_j:=\exp(-\beta
\delta_j), \ i=1,2, \ldots ,n-1, \ j=1,2, \ldots ,m-1$.
\end{thm}

\begin{rem}
\label{rem1} \
Observe that Fisher information on a single parameter ($\alpha$ or
$\beta$) depends only on the design points corresponding to that
particular parameter, e.g., $M_{\alpha}(n,m)=mM_{\alpha}(n)$, where
$M_{\alpha}(n)$ is the Fisher information corresponding to the
covariance parameter $\alpha$ of a stationary OU process
observed in design points $\{s_i, \ i=1,2,\ldots ,n\}$ of the interval
$[a_1,b_1]$.
\end{rem}

Now, with the help of  Theorem \ref{Mrn} one can formulate a result on
the restricted D-optimal design for the
parameters of the covariance structure of the OU sheet.

\begin{thm}
    \label{covpars} \
The restricted design which is D-optimal for estimation of the covariance
parameters $\alpha, \ \beta$ does not exist within the class of
admissible designs.
\end{thm}

From the point of view of a chemometrician, Theorem \ref{covpars} points out
that microergodicity should be added to the model in order to obtain
regular designs. Several ways are possible, for instance, nugget effect or
compounding \citep[see, e.g.,][]{MullerStehlik2}.

\begin{ex} \
  \label{ex:ex2.1}
Without loss of generality one may assume that the design space is
$\mathcal X=[0,1]^2$. Let $\alpha=0.6, \ \beta=1$, and consider the
case  $n=m=3$ where
$s_1=t_1=0, \ s_2:=d, \ t_2:=\delta, \ s_3=t_3=1$. For this particular
restricted
design we obviously have $d_1=d, \ d_2=1-d, \ \delta_1=\delta, \
\delta_2=1-\delta $. In Figure \ref{fig1}, where $\det \big(M_r(3,3)\big)$ is
plotted as function of $d$ and $\delta$, one can clearly see that the
maximal information is gained at the frontier points, when either
$d\in\{0,1\}$ or $\delta\in\{0,1\}$.
\begin{figure}[t]
\begin{center}
\leavevmode
\epsfig{file=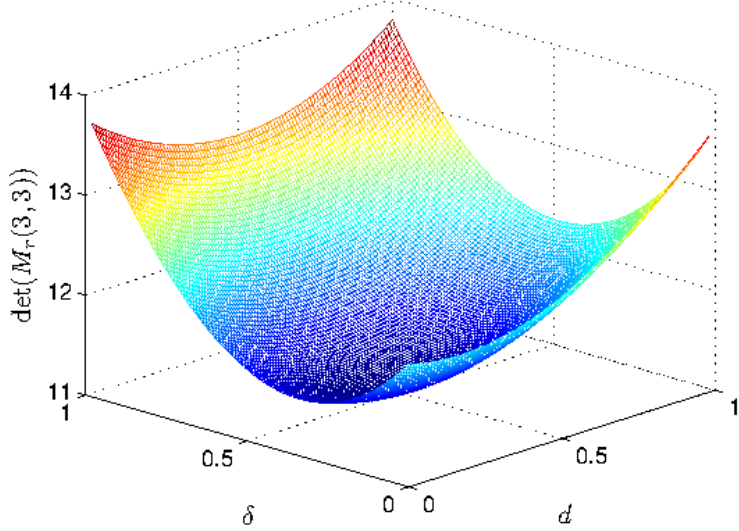,height=8cm}
\end{center}
\caption{Fisher information on correlation
  parameters $(\alpha,\beta)$ for $n=m=3$ as function of $d=d_1$ and
  $\delta=\delta_1$ in the case $\alpha=0.6, \ \beta=1$.}
\label{fig1}
\end{figure}
\end{ex}

Now, let us have a look at the free boundary directionally equidistant
designs, that is at designs where $d_1=d_2=\ldots =d_{n-1}=:d$ and
$\delta_1=\delta_2=\ldots =\delta_{m-1}=:\delta$. In this case a
D-optimal design is specified by directional distances $d$ and $\delta$ which
maximize
\begin{equation}
   \label{eq:eq4.3}
\det \big(M_r(n,m)\big)=\frac{(n-1)(m-1)d^2\delta^2}{\big({\mathrm
    e}^{2\alpha d}-1\big)^2\big({\mathrm
    e}^{2\beta\delta}-1\big)^2}\Big(nm \big({\mathrm
    e}^{2\alpha d}+1\big)\big({\mathrm e}^{2\beta \delta}+1\big)-4(n-1)(m-1)\Big).
\end{equation}
In the case of OU processes this
question does not appear, since for processes Fisher information on
covariance parameter based on $n$
equidistant design points depends linearly on the two-point design
Fisher information \citep{KS}.

\begin{thm}
  \label{freedes1}
If $nm\geq 2(n-1)(m-1)$ then  $\det \big(M_r(n,m)\big)$ is strictly
monotone decreasing both in $d$ and  $\delta$, so its maximum is
reached at $d=\delta=0$.  If  $nm<
2(n-1)(m-1)$ then for fixed and small enough $d$ ($\delta $), function $\det
\big(M_r(n,m)\big)$ has a single maximum in $\delta$ ($d$).
\end{thm}

\begin{figure}[t]
\begin{center}
\leavevmode
\epsfig{file=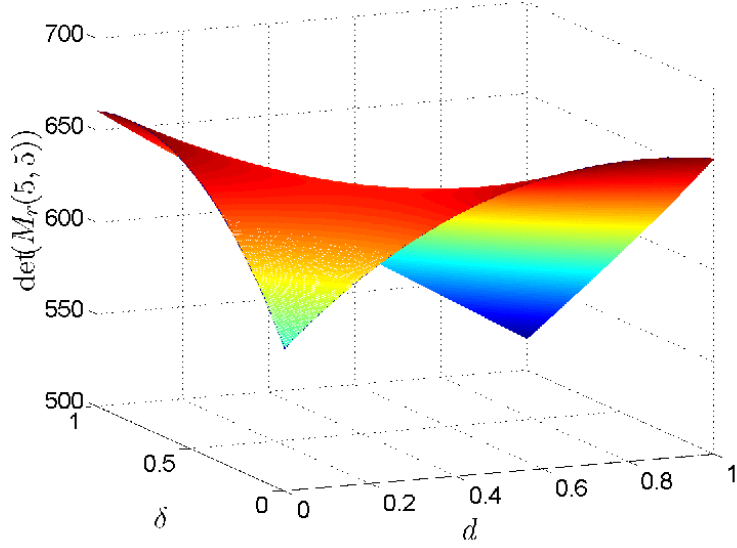,height=8cm}
\end{center}
\caption{Fisher information of boundary free design on correlation
  parameters $(\alpha,\beta)$ for $n=m=5$ in the case $\alpha=1, \ \beta=1$.}
\label{fig2}
\end{figure}

\begin{rem}
 Observe that for $1<n=m\in {\mathbb N}$ condition  $nm\geq 2(n-1)(m-1)$
 is equivalent to $n\leq 3$.
Further, if $nm\leq 2(n-1)(m-1)$ then the statement of Theorem
\ref{freedes1} does not imply the existence of a D-optimal
design. Figure \ref{fig2} shows that the extremal point of $\det
\big(M_r(n,m)\big)$ can be a saddle point and the maximum is reached
when either $d=0$ or $\delta=0$.
\end{rem}

\subsubsection{Estimation of all parameters}
Consider now the most general case, when both $\alpha,\ \beta$ and
$\theta$ are unknown and the Fisher information matrix on these
parameters equals
\begin{equation*}
M(n,m)=
\begin{bmatrix}
M_{\theta}(n,m) & 0 \\
0 &M_r(n,m)
\end{bmatrix},
\end{equation*}
where $M_{\theta}(n,m)$ and  $M_r(n,m)$  are Fisher information matrices on
$\theta$ and $r=(\alpha,\beta)^{\top}$, respectively, see
\eqref{eq:eq3.1} and \eqref{eq:eq4.1}. Thus, the objective function to
be maximized is $\det \big(M(n,m)\big)=M_{\theta}(n,m)\det
\big(M_r(n,m)\big)$.

\begin{ex} \
  \label{ex:ex2.2}
Consider the nine-point restricted design of Example \ref{ex:ex2.1},
that is  $\mathcal X=[0,1]^2$, \ $n=m=3$ and
$s_1=t_1=0, \ s_2:=d, \ t_2:=\delta, \ s_3=t_3=1$, implying $d_1=d, \
d_2=1-d, \ \delta_1=\delta, \ \delta_2=1-\delta $. In this case from
\eqref{eq:eq3.1} and \eqref{eq:eq4.2} we have
\begin{align}
\det &\,\big(M(3,3)\big)=\bigg(1+\frac{{\mathrm e}^{\alpha
    d}-1}{{\mathrm e}^{\alpha d}+1}+ \frac{{\mathrm e}^{\alpha
    (1-d)}-1}{{\mathrm e}^{\alpha (1-d)}+1}\bigg)
\bigg(1+\frac{{\mathrm e}^{\beta\delta}-1}{{\mathrm
    e}^{\beta\delta}+1}+ \frac{{\mathrm e}^{\beta(1-\delta)}
  -1}{{\mathrm e}^{\beta(1-\delta)}+1}\bigg) \nonumber \\
&\times
\Bigg(9\bigg(\frac{d^2\big({\mathrm e}^{2\alpha
    d}+1\big)}{\big({\mathrm e}^{2\alpha
    d}-1\big)^2}+\frac{(1-d)^2\big({\mathrm e}^{2\alpha
    (1-d)}+1\big)}{\big({\mathrm e}^{2\alpha (1-d)}-1\big)^2}  \bigg)
\bigg(\frac{\delta^2\big({\mathrm
    e}^{2\beta\delta}+1\big)}{\big({\mathrm
    e}^{2\beta\delta}-1\big)^2}+\frac{(1-\delta)^2\big({\mathrm
    e}^{2\beta(1-\delta)}+1\big)}{\big({\mathrm e}^{2\beta(1-\delta)}-1\big)^2}
\bigg)  \label{eq:eq4.4}\\
&-4\bigg(\frac{d}{{\mathrm e}^{2\alpha d}-1}+\frac{1-d}{{\mathrm
    e}^{2\alpha (1-d)}-1}\bigg)^2 \bigg(\frac{\delta}{{\mathrm
    e}^{2\beta\delta}-1}+\frac{1-\delta}{{\mathrm
    e}^{2\beta (1-\delta)}-1}\bigg)^2 \Bigg).\nonumber
\end{align}
Tedious calculations (see Section \ref{subs:subsA.5}) show  that
$\det\big(M(3,3)\big)$ 
has a single global minimum at $d=\delta=1/2$, while the maximum is
reached at the 
four vertices of $\mathcal X$, namely at $(0,0),\ (0,1), \ (1,0)$ and
$(1,1)$. In this way a restricted D-optimal design does not exist.
\end{ex}

Again, let us also have a look at the free boundary directionally equidistant
designs with directional distances $d$ and $\delta$. The objective
function to be maximized in order to get the D-optimal design is

\begin{align}
   \label{eq:eq4.5}
\det \big(M(n,m)\big)=&\frac{(n-1)(m-1)d^2\delta^2}{\big({\mathrm
    e}^{2\alpha d}-1\big)^2\big({\mathrm e}^{2\beta\delta}-1\big)^2\big({\mathrm
    e}^{\alpha d}+1\big)\big({\mathrm e}^{\beta \delta}+1\big)}
\big(n({\mathrm e}^{\alpha d}-1)+2\big)\big(m({\mathrm
    e}^{\beta\delta}-1)+2\big) \\
    &\times  \Big(nm \big({\mathrm
    e}^{2\alpha d}+1\big)\big({\mathrm e}^{2\beta
    \delta}+1\big)-4(n-1)(m-1)\Big). \nonumber
\end{align}
For simplicity assume $n=m$.

\begin{figure}[t]
\begin{center}
\leavevmode
\epsfig{file=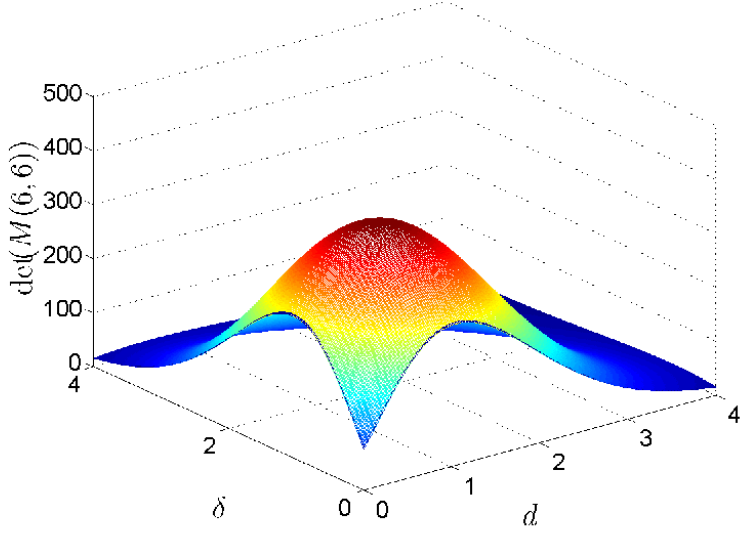,height=8cm}
\end{center}
\caption{Fisher information of boundary free design on all
  parameters  for $n=m=6$ in the case $\alpha=1, \ \beta=1$.}
\label{fig3}
\end{figure}

\begin{thm}
  \label{freedes2}
If $n=2$ then $\det \big(M(n,n)\big)$ is strictly monotone decreasing
both in $d$ and  $\delta$, so its maximum is reached at
$d=\delta=0$. If $n\geq 3 $ then  $\det \big(M(n,n)\big)$ has a global
maximum at $(d^*,\delta^*)$ which solves
\begin{equation}
   \label{eq:eq4.6}
n^2\big({\mathrm e}^{2\beta\delta}+1\big)g_1(\alpha
d,n)=4(n-1)^2g_2(\alpha d,n), \qquad
n^2\big({\mathrm e}^{2\alpha d}+1\big)g_1(\beta\delta
,n)=4(n-1)^2g_2(\beta\delta,n),  
\end{equation}
where
\begin{align}
g_1(x,n)&:={\mathrm e}^{5x}n(1-x)+{\mathrm e}^{4x}(2nx-3x-n+2)+
{\mathrm e}^{3x}x(1-4n)+{\mathrm e}^{2x}x(4n-7)+{\mathrm
  e}^x(x-n-nx)+n-2, \nonumber\\
g_2(x,n)&:={\mathrm e}^{3x}n(1-2x)+{\mathrm
  e}^{2x}(3nx-5x+2-n)+{\mathrm e}^x(x-n-nx)+n-2.   \label{eq:eq4.7}
\end{align}
\end{thm}

Theorem \ref{freedes2} shows that the situation here completely differs
from the case when only covariance parameters are estimated and an
optimal free boundary directionally equidistant design does exist. This
can clearly be observed on Figure \ref{fig3} showing $\det
\big(M(6,6)\big)$ for $\alpha=1, \ \beta=1$. Further, simulation
results show that for all $n\geq 3$ objective function $\det
\big(M(n,n)\big)$ has a unique maximal point (system \eqref{eq:eq4.6}
has a unique solution), however, a rigorous proof of this fact have
not been found yet.

\subsection{Optimal design with respect to IMSPE criterion}
 \label{subsec:subsec2.2}
As before, suppose we have observations $\big\{Y(s_i,t_j), \ i=1,2,\ldots,
n, \ j=1,2,\ldots, m\big\}$. The main aim of the kriging technique consists of
the prediction of the output
of the simulator on the experimental region. For any untried
location $(x_1,x_2)\in {\mathcal X}$ the estimation procedure is
focused on the best linear
unbiased estimator of $Y(x_1,x_2)$ given by $\widehat Y(x_1,x_2)=\widehat \theta
+R^{\top}(x_1,x_2)C^{-1}(n,m,r)({\mathbf Y}-{\mathbf 1}_{nm}\widehat
\theta)$, where ${\mathbf Y}=\big(Y(s_1,t_1), Y(s_1,t_2),\ldots
,Y(s_n,t_m)\big)^{\top}$ is the vector of observations, $\widehat\theta$ is the
generalized least squares estimator  of $\theta$, that is
$\widehat
\theta=\big({\mathbf 1}_{nm}^{\top}C^{-1}(n,m,r){\mathbf
  1}_{nm}\big)^{-1}{\mathbf 1}_{nm}^{\top} C^{-1}(n,m,r){\mathbf
  Y}$, and $R(x_1,x_2)$ is the vector of
correlations between $Y(x_1,x_2)$ and vector ${\mathbf Y}$ defined by
\!$R(x_1,x_2)\!\!=\!\!\big(\!\varrho(x_1,x_2,s_1,t_1),\ldots
,\varrho(x_1,x_2,s_i,t_j), \ldots
,\varrho(x_1,x_2,s_n,t_m)\!\big)\!^{\top}$\!\!,
where $\varrho(x_1,x_2,s_i,t_j):=\varrho_1(x_1,s_i)\varrho_2(x_2,t_j)$
with components
$\varrho_1(x_1,s_i):=\exp\big(-\alpha|x_1-s_i|\big)$ and
$\varrho_2(x_2,t_j):=\exp\big(-\beta|x_2-t_j|\big)$. Usually,
correlation parameters $\alpha, \beta$ are
unknown and will be estimated by maximum likelihood method. Thus,
the kriging predictor is obtained by substituting the maximum
likelihood estimators (MLE) $(\widehat
\alpha,\widehat \beta)$ for $(\alpha, \beta)$ and in such a case
$\widehat Y(x_1,x_2)$ is
called the MLE-empirical best linear unbiased predictor
\citep{Santner}.

In this way a natural criterion of optimality will  minimize suitable
functionals of the Mean Squared Prediction Error (MSPE) given by
\begin{equation}
   \label{eq:eq6.1}
\mspe\big(\widehat Y(x_1,x_2)\big):= \sigma^2
\Bigg[1-\big(1,\,R^{\top}(x_1,x_2)\big)
\left[
      \begin{BMAT}{c.c}{c.c}
      0 & {\mathbf 1}_{nm}^{\top}\\
      {\mathbf 1}_{nm} & C(n,m,r)
      \end{BMAT}
\right]^{-1}
\big(1,\,R^{\top}(x_1,x_2)\big)^{\top}\Bigg].
\end{equation}

Since the prediction accuracy is often
related to the entire prediction region $\mathcal{X}$ the
design criterion IMSPE is given by
\begin{equation*}
\imspe\big (\widehat
Y\big):=\sigma^{-2}\iint\limits_{\mathcal{X}}\mspe\big (\widehat
Y(x_1,x_2)\big)\,{\mathrm d}x_1\,{\mathrm d}x_2.
\end{equation*}

\begin{thm}
 \label{IMSPE}
\ Let us assume that the design space
$\mathcal{X}=[0,1]^2$ and since extrapolative prediction
is not advisable in kriging, we can set $s_1=t_1=0$ and $s_n=t_m=1$.
\begin{align}
\mspe\big (\widehat
Y(x_1,x_2)\big)=\sigma^2\Bigg[1-&\,\bigg(\varrho_1^2(x_1,s_n)+\sum_{i=1}^{n-1}
\frac{\big(\varrho_1(x_1,s_i)-\varrho_1(x_1,s_{i+1})p_i\big)^2}{1-p_i^2}\bigg)
\nonumber \\
&\times
\bigg(\varrho_2^2(x_2,t_m)+\sum_{j=1}^{m-1}
\frac{\big(\varrho_2(x_2,t_j)-\varrho_2(x_2,t_{j+1})q_j\big)^2}{1-q_j^2}\bigg)
\label{eq:eq6.2} \\
+\bigg(1+\sum_{i=1}^{n-1} \frac{1-p_i}{1+p_i}\bigg)^{-1}&
\bigg(1+\sum_{j=1}^{m-1} \frac{1-q_j}{1+q_j}\bigg)^{-1}
\Bigg(1-\bigg(\varrho_1(x_1,s_n)+\sum_{i=1}^{n-1}
\frac{\varrho_1(x_1,s_i)-\varrho_1(x_1,s_{i+1})p_i}{1+p_i}\bigg)
\nonumber \\
&\times
\bigg(\varrho_2(x_2,t_m)+\sum_{j=1}^{m-1}
\frac{\varrho_2(x_2,t_j)-\varrho_2(x_2,t_{j+1})q_j}{1+q_j}\bigg)\Bigg)^2\Bigg],
\nonumber
\end{align}
where again $p_i:=\exp(-\alpha d_i), \ q_j:=\exp(-\beta \delta _j)$
with $d_i:=s_{i+1}-s_i$ and
$\delta_j:=t_{j+1}-t_j, \ i=1,2, \ldots, n-1, \ j=1,2, \ldots ,m-1$.
Further,
\begin{align}
  \label{eq:eq6.3}
\imspe\big (\widehat Y\big)=&\,1-\bigg(\frac {n-1}{\alpha}-2\sum_{i=1}^{n-1}\frac
{d_ip_i^2}{1-p_i^2}\bigg)\bigg(\frac {m-1}{\beta}-2\sum_{j=1}^{m-1}\frac
{\delta_jq_j^2}{1-q_j^2}\bigg) \\
&+\bigg(1+\sum_{i=1}^{n-1}
\frac{1-p_i}{1+p_i}\bigg)^{-1}\bigg(1+\sum_{j=1}^{m-1}
\frac{1-q_j}{1+q_j}\bigg)^{-1}\Bigg[1-\frac
8{\alpha\beta}\bigg(\sum_{i=1}^{n-1}\frac
{1-p_i}{1+p_i}\bigg)\bigg(\sum_{j=1}^{m-1}\frac
{1-q_j}{1+q_j}\bigg) \nonumber \\
&+\bigg(\sum_{i=1}^{n-1}\frac
{1-p_i^2+2\alpha d_ip_i}{\alpha (1+p_i)^2}\bigg)\bigg(\sum_{j=1}^{m-1}\frac
{1-q_j^2+2\beta \delta_jq_j}{\beta (1+q_j)^2}\bigg)\Bigg]. \nonumber
\end{align}
For any sample size the  directionally equidistant design
$d_1=d_2=\ldots =d_{n-1}$ and $\delta_1=\delta _2=\ldots =\delta_{m-1}$  is
optimal with respect to the IMSPE criterion.

\end{thm}

\begin{rem} \ We remark that \eqref{eq:eq6.3} is an extension
  of the IMSPE criterion for the classical OU process given by
  \citet[Proposition 4.1]{baz}, while the optimality result
  generalizes Proposition 4.2 of \citet{baz}.
\end{rem}

\subsection{Optimal design with respect to entropy criterion}
  \label{subsec:subsec2.3}
Another possible approach to optimal design is to find locations which maximize
the amount of obtained information. Following the ideas of \citet{sw}
one has to maximize the entropy $\ent ({\mathbf Y})$ of the
observations corresponding to the chosen design, which in the Gaussian
case form an $nm$-dimensional normal vector with covariance
matrix $\sigma ^2\,C(n,m,r)$, that is
\begin{equation*}
\ent ({\mathbf Y})=\frac {nm}2\big(1+\ln (2\pi \sigma^2)\big)+\frac
12\ln\det C(n,m,r).
\end{equation*}

\begin{thm}
  \label{entropy}
In our setup entropy $\ent ({\mathbf Y})$
has the form
\begin{equation}
  \label{eq:eq7.1}
\ent ({\mathbf Y})=\frac {nm}2\big(1+\ln (2\pi \sigma^2)\big)+\frac
m2 \sum_{i=1}^{n-1} \ln \big(1-p_i^2\big)+\frac
n2 \sum_{j=1}^{m-1} \ln \big(1-q_j^2\big).
\end{equation}
For any sample size the  directionally equidistant design
$d_1=d_2=\ldots =d_{n-1}$ and $\delta_1=\delta _2=\ldots =\delta_{m-1}$
is optimal with respect to the entropy criterion.
\end{thm}

\section{D-optimal designs for the Arrhenius model with OU error}
  \label{sec:sec3}
In the present section we derive objective functions for D-optimal
designs for estimating parameters of the Arrhenius model \eqref{eq:GEmodel}.
We consider the stationary process
\begin{equation}
   \label{Arhmodel}
Y(s,t) =\big( A/t^{\mu}\big){\mathrm e}^{-B/t} +\varepsilon (s,t),
\end{equation}
observed on a compact design space
$\mathcal{X}=[a_1,b_1]\times [a_2, b_2]$, where $b_1>a_1$ and
$b_2>a_2$ and $\varepsilon (s,t), \ s,t\in {\mathbb R}$, is again a
stationary Ornstein-Uhlenbeck
sheet, that is a zero mean Gaussian process with covariance structure
(\ref{oucovmod}). Since parameter $A$ is usually known, without loss
of generality we may assume $A=1$ and consider model \eqref{Arhmodel}
with trend function $\eta(s,t;\mu,B):=\big(1/t^{\mu}\big){\mathrm
  e}^{-B/t}$.

From the point of view of applications we distinguish two important cases.
\begin{itemize}
\item Rate $\mu$ is known, which is an assumption made by several
  authors, see, e.g., \citet{Heb87}. The uncorrelated case has already been
  studied by \citet{RodSan09}, where the authors proved that for
  approximated designs a two-point design is optimal.
\item Rate $\mu$ is unknown and one has to estimate it together with
  $B$. For this model the uncorrelated case has also been studied,
  \citet{Chemo12} considered both equidistant and general designs.
\end{itemize}

\subsection{Estimation of trend}
  \label{subs:subs3.1}
Assume that covariance parameters $\alpha, \beta$ and $\sigma$ of the
OU sheet and rate $\mu$ of the Arrhenius model are given and we are
interested in estimation of the trend parameter $B$. The Fisher
information on $B$ based on observations $\big\{Y(s_i,t_j), \ i=1,2,\ldots,
n, \ j=1,2,\ldots, m\big\}$ of the process \eqref{Arhmodel} equals
$M_{B}(n,m)=F^{\top}(n,m,B)C^{-1}(n,m,r)F^{\top}(n,m,B)$, where
\begin{equation*}
F(n,m,B):=\bigg (\frac {\eta(s_1,t_1;\mu,B)}{\partial B}, \frac
{\eta(s_1,t_2;\mu,B)}{\partial B}, \ldots, \frac
{\eta(s_n,t_m;\mu,B)}{\partial B}\bigg)^{\top}.
\end{equation*}

\begin{thm}
  \label{arhfish}
In our setup
\begin{equation}
  \label{arhfishform}
M_{B}(n,m)=\bigg(1+\sum_{i=1}^{n-1}\frac{1-p_i}{1+p_i}\bigg)\bigg(\kappa_m^2
+\sum_{j=1}^{m-1}\frac{(\kappa_j-\kappa_{j+1}q_j)^2}{1-q_j^2}\bigg),
\end{equation}
where $\kappa_j:=-\exp\big(-B/t_j\big)/t_j^{\mu+1}$ if $t_j\ne 0$, and
$\kappa_j:=0$, otherwise.
\end{thm}

In case one has to estimate both $\mu$ and $B$, the objective function
to be maximized in order to get the D-optimal design is $\det
\big(M_{\mu,B}(n,m)\big)$, where again
$M_{\mu, B}(n,m)=G^{\top}(n,m,\mu,B)C^{-1}(n,m,r)G^{\top}(n,m,\mu,B)$ with
\begin{equation*}
G(n,m,\mu,B):=\begin{bmatrix}
\frac {\eta(s_1,t_1;\mu,B)}{\partial \mu}& \frac
{\eta(s_1,t_2;\mu,B)}{\partial \mu}& \ldots & \frac
{\eta(s_n,t_m;\mu,B)}{\partial \mu}\\
\frac {\eta(s_1,t_1;\mu,B)}{\partial B}& \frac
{\eta(s_1,t_2;\mu,B)}{\partial B}& \ldots & \frac
{\eta(s_n,t_m;\mu,B)}{\partial B}
\end{bmatrix}^{\top}.
\end{equation*}
\begin{thm}
  \label{arhfishmu}
In our setup
\begin{align}
  \label{arhfishmuform}
M_{\mu,B}(n,m)=&\bigg(1+\sum_{i=1}^{n-1}\frac{1-p_i}{1+p_i}\bigg) \\
&\times
\begin{bmatrix}
\lambda_m^2
+\sum_{j=1}^{m-1}\frac{(\lambda_j-\lambda_{j+1}q_j)^2}{1-q_j^2}&
\lambda_m\kappa_m
+\sum_{j=1}^{m-1}\frac{(\lambda_j-\lambda_{j+1}q_j)
  (\kappa_j-\kappa_{j+1}q_j)}{1-q_j^2}\\
\lambda_m\kappa_m
+\sum_{j=1}^{m-1}\frac{(\lambda_j-\lambda_{j+1}q_j)
  (\kappa_j-\kappa_{j+1}q_j)}{1-q_j^2}&
\kappa_m^2
+\sum_{j=1}^{m-1}\frac{(\kappa_j-\kappa_{j+1}q_j)^2}{1-q_j^2}
\end{bmatrix}, \nonumber
\end{align}
where $\kappa_j$ is the same quantity as in Theorem \ref{arhfish}, while
$\lambda_j:=-\log(t_j)\exp\big(-B/t_j\big)/t_j^{\mu}$ if $t_j\ne 0$, and
$\lambda_j:=0$, otherwise.
\end{thm}

Theorems \ref{arhfish} and \ref{arhfishmu} show that for estimating merely
the trend parameters one can treat the two coordinate directions
separately. Hence, in the first coordinate direction the
maximum is reached with the equidistant design
$d_1=d_2=\ldots=d_{n-1}$, while in the second
direction one can consider, e.g., the results of  \citet{Chemo12} for
the classical OU process.

\begin{ex}
  \label{ex:ex3.1}
Consider  a four point grid design, i.e. $n=m=2$. Without loss of generality
we may assume
$s_1=t_1=0$ implying $s_2=d$ and $t_2=\delta$.
In this case the Fisher information \eqref{arhfishform} on $B$ equals
\begin{equation*}
M_B(2,2)=\frac 2{1-\exp(-\alpha d)}
\frac{\exp(-2B/\delta)}{\big(1-\exp(-2\beta \delta)\big)\delta^{2(\mu+1)}},
\end{equation*}
which function is monotone increasing in its first variable $d$.
Further, short calculation shows that if $\mu>-1$ then the maximum in
$\delta$ is attained at the unique solution of the equation
\begin{equation*}
\big(B-(\mu+1)\delta\big)\big(\exp(2\beta\delta)-1\big)=\beta\delta^2.
\end{equation*}

In case $\mu<-1$, that is in particular interesting for chemometricians,
one can employ the maximin approach  \citep[see, e.g.,][]{MAXIMIN} which
seeks designs maximizing the minimum of the design criterion. In
our case this means maximization of
\begin{equation}
  \label{MAXIMINcriteria}
\min_{\alpha,\beta>0}
M_B(2,2)=2\exp(-2B\delta)\delta^{-2(\mu+1)}.
\end{equation}
Obviously, if $\mu<-1$ then the maximum of \eqref{MAXIMINcriteria} is
reached at $\delta^*=-(\mu+1)/B$. Although the maximization of
\eqref{MAXIMINcriteria} is pretty easy, one should take care about the
interpretation of such a result as, e.g., the optimal design does not
depend on $d$.

Maximin approach, anyhow, cannot be automatized without further
considerations since, for instance, maximin designs are of no relevance for
criteria, where design distances are multiplied  by some nuisance
parameters, see, e.g., \eqref{eq:eq3.1}.
\end{ex}

\begin{rem} Under the conditions of Example
  \ref{ex:ex3.1} ($s_1=t_1=0$) we have
  $\det\big(M_{\mu,B}(2,2)\big)=0$, that is the four point grid design
  does not provide information on trend parameters $\mu$ and $B$.
\end{rem}

\subsection{Estimation of all parameters}
  \label{subs:subs3.2}
Assume first that the rate $\mu$ is known and one has to
estimate trend parameter $B$ and covariance parameters
$(\alpha,\beta)$. Obviously, the Fisher information matrix on these
parameters  based on observations $\big\{Y(s_i,t_j), \ i=1,2,\ldots,
n, \ j=1,2,\ldots, m\big\}$ of the process \eqref{Arhmodel} equals
\begin{equation*}
{\mathcal M}(n,m)=
\begin{bmatrix}
M_B(n,m) & 0 \\
0 &M_r(n,m)
\end{bmatrix},
\end{equation*}
where $M_B(n,m)$ and  $M_r(n,m)$  are defined by \eqref{arhfishform}
and \eqref{eq:eq4.1}, respectively. Hence, in order to obtain a
D-optimal design one has to maximize $\det\big({\mathcal
  M}(n,m)\big)=M_B(n,m) \det\big(M_r(n,m)\big)$.

\begin{ex}
 \label{ex:ex3.2}
Consider again the settings of Example \ref{ex:ex3.1}, that is a four
point grid design ($n=m=2$) under the assumption $s_1=t_1=0$. In this
case we have
\begin{equation*}
{\mathcal M}(2,2)=
\frac{8d^2\exp(-2B/\delta)\exp(-2\beta\delta)\exp(-2\alpha
  d)\big(1+\exp(-2\alpha d)+\exp(-2\beta
  \delta)\big)}{\delta^{2\mu}\big(1-\exp(-2\beta
  \delta)\big)^3\big( 1-\exp(-2\alpha d)\big)^2
  \big(1+\exp(-\alpha d)\big)}, \qquad  d,\delta \geq 0.
\end{equation*}
Tedious calculations (see Section \ref{subs:subsA.11}) show  that for
$d,\delta\geq0$
function ${\mathcal M}(2,2)$ is monotone decreasing
in $d$, while in $\delta$ it has a maximum at the unique solution of
the equation
\begin{equation*}
\beta \delta^2-\mu\delta +B+{\mathrm
  e}^{2\beta\delta}\big(2\beta(2+p^2)\delta^2+(B-\mu\delta)p^2 \big)+{\mathrm
  e}^{4\beta\delta}(1+p^2)\big(\beta \delta^2 +\mu\delta -B\big)=0.
\end{equation*}
Hence, the optimal four point grid design collapses in its first
coordinate.
\end{ex}

If rate $\mu$ is also unknown, the Fisher information matrix on
$(\mu,B,\alpha,\beta)$ based on $\big\{Y(s_i,t_j), \ i=1,2,\ldots,
n, \ j=1,2,\ldots, m\big\}$ equals
\begin{equation*}
{\mathfrak M}(n,m)=
\begin{bmatrix}
M_{\mu,B}(n,m) & 0 \\
0 &M_r(n,m)
\end{bmatrix},
\end{equation*}
where $M_{\mu,B}(n,m)$ and  $M_r(n,m)$  are defined by \eqref{arhfishmuform}
and \eqref{eq:eq4.1}, respectively. In this case the D-optimal design
maximizes objective function $\det\big({\mathfrak
  M}(n,m)\big)=\det\big(M_{\mu,B}(n,m)\big) \det\big(M_r(n,m)\big)$.

\section{Comparisons of designs}
   \label{sec:sec4}
Methane emissions compose a very complicated process which mixes
stochasticity with chaos \citep[see, e.g.,][]{Addiscott,Sabol}, thus fitting of
two dimensional OU sheet could be a remedy to several problems
which occurred in univariate settings \citep{Chemo12}.
In this section we provide efficiency comparisons for selected
important methane kinetic reactions, both in standard (Earth) and
non-standard (troposphere) conditions.
The current work is the first comprehensive comparison of
  the statistical information of designs for OU sheets, which gives
  its novelty both methodologically and from the point of view of
  applications.
\begin{table}[t!]
\begin{center}\footnotesize
\begin{tabular}{|c|c|c|c|c|c|}
\hline
\multicolumn{2}{|c|}{}&$\alpha=0.001, \ \beta=0.01$&$\alpha=0.1, \
\beta=1$&$\alpha=1, \ \beta=1$&$\alpha=1, \ \beta=10$ \\\hline
&monotonic&1.3118&29.8651&61.2545&63.9937\\
$D-opt.$&rectangular&1.3328&57.4388&63.7483&64.00\\
&rel. eff. (\%)&98.43&51.99&96.09&99.99\\\hline
&monotonic&-33.0446&86.1318&90.7964&90.8121\\
$entropy$&rectangular&-51.1507&90.7111&90.8119&90.8121\\
&rel. eff. (\%)&64.60&94.95&99.98&100\\\hline
\end{tabular}
\end{center}
\caption{$M_{\theta} (n,m)$ and entropy values corresponding to the
  optimal monotonic and to the
  rectangular grid design and relative efficiency of the optimal
  monotonic design.}
\label{tab1}
\end{table}

\subsection{Comparisons of  designs for tropospheric methane
  measurements}
  \label{subs:subs4.1}

As discussed by \citet{Lelieveld}, tropospheric methane measurements
are fundamental for climate change models and \citet{Vag91} utilized a
$62$ point design to measure the tropospheric methane flux.
In Theorem \ref{trend} the exact form of $M_{\theta} (n,m)$ is derived
only for restricted regular designs, however, one might ask what is the
relative efficiency of the optimal value of $M_{\theta} (n)$ on
monotonic sets \citep{bm} containing $n\times m$ design points
compared to the $M_{\theta} (n,m)$ of a rectangular
grid with the same number of points. Since the designs for methane
used in \citet{Vag91} typically have around 62 points,
we should consider a $64$ point design comparison of, e.g., a $8\times 8$
regular grid with a $64$ points monotonic set for covariance
parameters  $\alpha,\beta \in \{0.001,0.01,0.1,1,10\}$ and
design space  $[223,420]\times [0.84,43.51]$.

\begin{figure}[t]
\begin{center}
\leavevmode
\epsfig{file=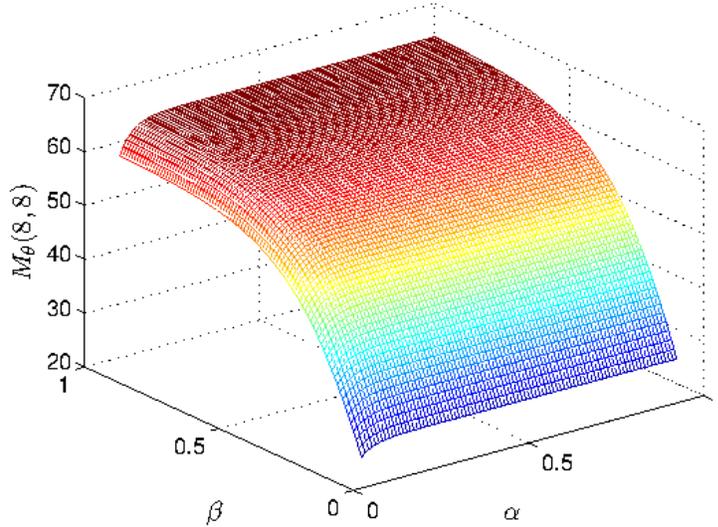,height=8cm}
\end{center}
\caption{Fisher information on $\theta$ as functions of correlation
  parameters $(\alpha,\beta)$ for $n=8$ and $m=8$.}
\label{fig4}
\end{figure}

Table \ref{tab1} gives
the optimal values of $M_{\theta}(64)$ on monotonic sets,
$M_{\theta}(8,8)$ values for regular designs and the relative
efficiencies of the optimal $M_{\theta}(64)$ values
on monotonic sets for different combinations of
parameters $(\alpha,\beta)$. Observe, that for $\alpha=0.1, \,
\beta=1$ the optimal monotonic design gives much lower values of
Fisher information on $\theta$ than the regular grid, while for
the other combinations of parameters the relative efficiency is
slightly below $100\,\%$. For the entropy criterion we obtain the same
results. In Figure \ref{fig4} the optimal value of Fisher
information on $\theta$ is plotted as a function of correlation
parameters $(\alpha,\beta)$ for $n=8$ and $m=8$.

\subsection{Comparisons of  designs for the rate of methane reactions
  with $OH$}
  \label{subs:subs4.2}
The growth rate of tropospheric methane is determined by the balance
between surface emissions and photo-chemical destruction by the
hydroxyl radical $OH,$ the major atmospheric oxidant. Such reaction
can happen at various temperature modes, for instance,
\citet{Bon02} measured the rate constants of the reactions of $OH$
radicals with methane in the temperature range $295-618K$. The
following 4 tables provide efficiency of original 
designs of \citet{Bon02} together with efficiencies of monotonic and 
regular grid designs $3\times 3, \ 2\times 5, \ 5\times 2, \ 3\times
4, \ 4\times 3, \ 3\times 2$ and $2\times 3$, respectively. Tables
\ref{tab2}--\ref{tab5} utilize the setups described in Tables 1--4 of
\citet{Bon02}. As one can see, in most of the situations monotonic and
regular grid designs outperform the original designs.

\begin{table}[b!]
\begin{center}\footnotesize
\begin{tabular}{|c|c|c|c|c|c|c|}
\hline
\multicolumn{2}{|c|}{}&$\alpha\!=\!0.001,
\beta\!=\!0.01$&$\alpha\!=\!0.1, \beta\!=\!0.01$&$\alpha\!=\!0.1,
\beta\!=\!1$&$\alpha\!=\!1, \beta\!=\!1$&$\alpha\!=\!1, \beta\!=\!10$
\\\hline
&\citet{Bon02}&3.1261&8.7785&8.9904&9.0000&9.0000\\
$D-opt.$&mon., $n=9$&3.2067&8.9107&9.0000&9.0000&9.0000\\
&$3\times 3$ r.grid&3.0305&7.6660&9.0000&9.0000&9.0000\\\hline
&\citet{Bon02}&9.8567&12.7665&12.7704&12.7704&12.7704\\
$Ent.$&mon., $n=9$&11.2150&12.7703&12.7704&12.7704&12.7704\\
&$3\times 3$ r.grid&9.2225&12.7231&12.7704&12.7704&12.7704\\\hline
\end{tabular}

\end{center}
\caption {$M_\theta(n,m)$ and entropy values corresponding to the
  optimal monotonic and to the rectangular grid design  together with values of
optimality criteria for measurements given in \citet[Table 1]{Bon02}.}
\label{tab2}
\end{table}

\begin{table}[t!]
\begin{center}\footnotesize
\begin{tabular}{|c|c|c|c|c|c|c|}
\hline
\multicolumn{2}{|c|}{}&$\alpha\!=\!0.001,
\beta\!=\!0.01$&$\alpha\!=\!0.1, \beta\!=\!0.01$&$\alpha\!=\!0.1,
\beta\!=\!1$&$\alpha\!=\!1, \beta\!=\!1$&$\alpha\!=\!1, \beta\!=\!10$
\\\hline
&\citet{Bon02}&1.1853&6.9087&7.0855&8.7813&9.2477\\
$D-opt.$&mon., $n=10$&1.1858&9.5186&9.7151&10.0000&10.0000\\
&$2\times 5$ r.grid&1.1884&2.0487&6.3460&6.3460&9.9999\\
&$5\times 2$ r.grid&1.1897&5.1192&9.9189&9.9239&10.0000\\\hline
&\citet{Bon02}&-0.8169&11.9268&12.5103&14.0660&14.1336\\
$Ent.$&mon., $n=10$&2.7830&14.1860&14.1882&14.1894&14.1894\\
&$2\times 5$ r.grid&-2.5767&-0.7201&13.8227&13.8227&14.1894\\
&$5\times 2$ r.grid&0.6346&8.2463&14.1892&14.1892&14.1894\\\hline
\end{tabular}

\end{center}
\caption {$M_\theta(n,m)$ and entropy values corresponding to the
  optimal monotonic and to the rectangular grid design  together with values of
optimality criteria for measurements given in \citet[Table 2]{Bon02}.}
\label{tab3}
\end{table}

\begin{table}[t!]
\begin{center}\footnotesize
\begin{tabular}{|c|c|c|c|c|c|c|}
\hline
\multicolumn{2}{|c|}{}&$\alpha\!=\!0.001,
\beta\!=\!0.01$&$\alpha\!=\!0.1, \beta\!=\!0.01$&$\alpha\!=\!0.1,
\beta\!=\!1$&$\alpha\!=\!1, \beta\!=\!1$&$\alpha\!=\!1, \beta\!=\!10$
\\\hline
&\citet{Bon02}&1.1816&6.7348&6.9218&7.6265&9.0242\\
$D-opt.$&mon., $n=12$&1.1818&10.8570&11.2215&12.0000&12.0000\\
&$3\times 4$ r.grid&1.1850&3.0669&8.6804&8.6804&12.0000\\
&$4\times 3$ r.grid&1.1852&4.0890&10.4462&10.4466&12.0000\\\hline
&\citet{Bon02}&-5.7821&3.0845&12.3312&12.9532&16.4642\\
$Ent.$&mon., $n=12$&1.9060&17.0107&17.0199&17.0273&17.0273\\
&$3\times 4$ r.grid&-4.0505&1.1408&16.7911&16.7911&17.0273\\
&$4\times 3$ r.grid&-2.9378&4.4983&16.9807&16.9807&17.0273\\\hline
\end{tabular}

\end{center}
\caption {$M_\theta(n,m)$ and entropy values corresponding to the
  optimal monotonic and to the rectangular grid design together with values of
optimality criteria for measurements given in \citet[Table 3]{Bon02}.}
\label{tab4}
\end{table}

\begin{table}[t!]
\begin{center}\footnotesize
\begin{tabular}{|c|c|c|c|c|c|c|}
\hline
\multicolumn{2}{|c|}{}&$\alpha\!=\!0.001,
\beta\!=\!0.01$&$\alpha\!=\!0.1, \beta\!=\!0.01$&$\alpha\!=\!0.1,
\beta\!=\!1$&$\alpha\!=\!1, \beta\!=\!1$&$\alpha\!=\!1, \beta\!=\!10$
\\\hline
&\citet{Bon02}, $n=7$&1.0057&1.1531&1.5630&2.2240&4.5042\\
&\citet{Bon02}, $n=6$&1.0057&1.1531&1.5630&2.2240&4.4850\\
$D-opt.$&mon.,$n=7$&1.0057&1.1542&1.5683&2.8570&5.4387\\
&mon.,$n=6$&1.0057&1.1542&1.5675&2.8309&5.0721\\
&$2\times 3$ r.grid&1.0057&1.1537&1.6244&2.6938&5.6029\\
&$3\times 2$ r.grid&1.0057&1.1545&1.6061&3.1714&4.5396\\\hline
&\citet{Bon02}, $n=7$&-8.3075&-6.4357&4.9754&5.1821&8.9398\\
&\citet{Bon02}, $n=6$&-5.4333&-3.5616&5.5473&5.7539&8.3806\\
$Ent.$&mon., $n=7$&-6.7914&2.9548&6.4778&8.9552&9.8647\\
&mon., $n=6$&-4.9681&3.1294&6.0021&7.9077&8.4873\\
&$2\times 3$ r.grid&-8.7323&-2.2476&6.1896&7.3797&8.5095\\
&$3\times 2$ r.grid&-9.2498&-0.3290&5.5038&8.1021&8.4115\\\hline
\end{tabular}

\end{center}
\caption {$M_\theta(n,m)$ and entropy values corresponding to the
  optimal monotonic and to the rectangular grid design together with values of
optimality criteria for measurements given in \citet[Table 4]{Bon02}.}
\label{tab5}
\end{table}

\citet{Dun93} measured absolute rate coefficients for the reactions of $OH$
radical with $CH_4$ ($k_1$) and perdeuterated methane $d_4$ ($k_2$.)
Authors characterized $k_1$ and $k_2$ over the temperature range $293-800K.$
Finally, they found an excellent agreement of their results with
determinations of $k_1$ at lower temperatures of \citet{Vag91}. Now,
let us consider rates $k_1$ and $k_2$ of Table 1 of \citet{Dun93}.
We obtain the following comparisons (Table \ref{tab6}-\ref{tab7}) of
efficiencies of the monotonic and $2\times 5$ and $5\times 2 $ regular grid
designs with the original designs of \citet{Dun93}.
These results show that in most of the cases, the monotonic and regular grid
designs are more efficient than the original one.

\begin{table}
\begin{center}\footnotesize
\begin{tabular}{|c|c|c|c|c|c|c|}
\hline
\multicolumn{2}{|c|}{}&$\alpha\!=\!0.001,
\beta\!=\!0.01$&$\alpha\!=\!0.1, \beta\!=\!0.01$&$\alpha\!=\!0.1,
\beta\!=\!1$&$\alpha\!=\!1, \beta\!=\!1$&$\alpha\!=\!1, \beta\!=\!10$
\\\hline
&\citet{Dun93}&4.5728&9.4857&9.9959&10.0000&10.0000\\
$D-opt.$&mon., $n=10$&4.7604&9.9721&10.0000&10.0000&10.0000\\
&$2\times 5$ r.grid&4.9144&7.8743&10.0000&10.0000&10.0000\\
&$5\times 2$ r.grid&2.5049&9.9944&9.9999&10.0000&10.0000\\\hline
&\citet{Dun93}&12.2328&14.1366&14.1894&14.1894&14.1894\\
$Ent.$&mon., $n=10$&13.3584&14.1894&14.1894&14.1894&14.1894\\
&$2\times 5$ r.grid&12.9678&14.0944&14.1894&14.1894&14.1894\\
&$5\times 2$ r.grid&8.2035&14.1894&14.1894&14.1894&14.1894\\\hline
\end{tabular}

\end{center}
\caption {$M_\theta(n,m)$ and entropy values corresponding to the
optimal monotonic and to the rectangular grid design together with values of
optimality criteria for $k_1$ measurements given in \citet[Table 1]{Dun93}.}
\label{tab6}
\end{table}

\begin{table}
\begin{center}\footnotesize
\begin{tabular}{|c|c|c|c|c|c|c|}
\hline
\multicolumn{2}{|c|}{}&$\alpha\!=\!0.001,
\beta\!=\!0.01$&$\alpha\!=\!0.1, \beta\!=\!0.01$&$\alpha\!=\!0.1,
\beta\!=\!1$&$\alpha\!=\!1, \beta\!=\!1$&$\alpha\!=\!1, \beta\!=\!10$
\\\hline
&\citet{Dun93}&3.0778&11.7720&11.9798&12.0000&12.0000\\
$D-opt.$&mon., $n=12$&3.1465&11.8465&12.0000&12.0000&12.0000\\
&$3\times 4$ r.grid&3.3749&8.0858&12.0000&12.0000&12.0000\\
&$4\times 3$ r.grid&3.1184&9.9557&12.0000&12.0000&12.0000\\\hline
&\citet{Dun93}&11.2608&17.0260&17.0272&17.0273&17.0273\\
$Ent.$&mon., $n=12$&13.7036&17.0270&17.0273&17.0273&17.0273\\
&$3\times 4$ r.grid&12.9774&16.6656&17.0273&17.0273&17.0273\\
&$4\times 3$ r.grid&11.3202&16.9405&17.0273&17.0273&17.0273\\\hline
\end{tabular}

\end{center}
\caption {$M_\theta(n,m)$ and entropy values corresponding to the
optimal monotonic and to the rectangular grid design together with values of
optimality criterion for $k_2$ measurements given in \citet[Table 2]{Dun93}.}
\label{tab7}
\end{table}

\section{Conclusions}

\label{sec:sec5}

Both Kyoto protocol \citep{Lelieveld} and recent Scandinavian and
Polish summits
in 2013 pointed out necessity to develop precise statistical modelling
of climate change. This, in particular should be addressed by
developing of optimal, or at least benchmarking designs for 
complex climatic models. The current work aims to contribute here for the
case of methane modelling in troposphere, lowest part of atmosphere.
As can be well seen in the paper, optimal designs for univariate case
(OU process, see \citet{Chemo12}) and planar OU sheets differ. Obviously,
planar OU sheet is much more precise, since it allows
variability both in temperature (main chemically understood driver of
chemical kinetics) and in a second variable, which can be either atmospheric
pressure or any other relevant quantity. Temperature itself is also
regressor, i.e. variable entering into trend parameter $k_1$.
One valuable further research direction, enabled by the second variable ``$s$''
will be direct modelling of reaction kinetics.
The optimal design for spatial process of methane flux  can be helpful
for better understanding the emerging issues of paleoclimatology
\citep{McShane}, which in major part relates to large variability.

\bigskip

\noindent {\bf Acknowledgment.} This research has been supported by
the Hungarian --Austrian
intergovernmental S\&T cooperation program T\'ET\_{}10-1-2011-0712
and partially supported the T\'AMOP-4.2.2.C-11/1/KONV-2012-0001
project. The project has been supported by the European Union,
co-financed by the European Social Fund. M. Stehl\'\i k acknowledges the
support of ANR project Desire FWF I 833-N18
and Fondecyt Proyecto Regular N° 1151441.  K. Sikolya has been supported
by T\'AMOP 4.2.4. A/2-11-1-2012-0001 project ``National Excellence
Program -- Elaborating and operating an inland student and researcher
personal support system''. The project was subsidized by the European
Union and co-financed by the European Social Fund.

\begin{appendix}
\section{Appendix}
  \label{sec:secA}
  
\subsection{Proof of Theorem \ref{trend}}
      \label{subs:subsA.1}

According to the notations of Section 
\ref{subsec:subsec2.1} let
$d_i:=s_{i+1}-s_i, \ \delta _j:=t_{j+1}-t_j$ and $p_i:=\exp(-\alpha
d_i), $ \\$\ q_j:=\exp(-\beta \delta _j)$. Short calculation shows that
\begin{equation}
   \label{eq:eqA.1}
C(n,m,r)=P(n,r)\otimes Q(m,r),
\end{equation}
where
\begin{align*}
P(n,r):=&\,
     \begin{bmatrix}
        1 &p_1 &p_1p_2  &p_1p_2p_3 &\dots &\dots &\prod_{i=1}^{n-1}p_i \\
        p_1 &1 &p_2 &p_2p_3 &\dots &\dots &\prod_{n=2}^{n-1}p_i \\
        p_1p_2 &p_2 &1 &p_3 &\dots &\dots &\prod_{i=3}^{n-1}p_i \\
        p_1p_2p_3 &p_2p_3 &p_3 &1 &\dots &\dots &\vdots \\
        \vdots &\vdots &\vdots &\vdots &\ddots & &\vdots \\
        \vdots &\vdots &\vdots &\vdots & &\ddots &p_{n-1} \\
        \prod_{i=1}^{n-1}p_i &\prod_{i=2}^{n-1}p_i
        &\prod_{i=3}^{n-1}p_i &\dots &\dots &p_{n-1} &1\\
     \end{bmatrix},\\[2mm]
Q(m,r):=&\,
     \begin{bmatrix}
        1 &q_1 &q_1q_2  &q_1q_2q_3 &\dots &\dots &\prod_{j=1}^{m-1}q_j \\
        q_1 &1 &q_2 &q_2q_3 &\dots &\dots &\prod_{j=2}^{m-1}q_j \\
        q_1q_2 &q_2 &1 &q_3 &\dots &\dots &\prod_{j=3}^{m-1}q_j \\
        q_1q_2q_3 &q_2q_3 &q_3 &1 &\dots &\dots &\vdots \\
        \vdots &\vdots &\vdots &\vdots &\ddots & &\vdots \\
        \vdots &\vdots &\vdots &\vdots & &\ddots &q_{m-1} \\
        \prod_{j=1}^{m-1}q_j &\prod_{j=2}^{m-1}q_j
        &\prod_{j=3}^{m-1}q_j &\dots &\dots &q_{m-1} &1\\
     \end{bmatrix}.
\end{align*}
By the properties of the Kronecker product
\begin{equation}
   \label{eq:eqA.2}
C^{-1}(n,m,r)=P^{-1}(n,r)\otimes Q^{-1}(m,r),
\end{equation}
and, according to the results of \citet{KS}, e.g., the inverse of
$P(n,r)$ equals
\begin{equation}
   \label{eq:eqA.3}
  P^{-1}(n,r)=
     \begin{bmatrix}
        \frac{1}{1-p_1^2} &\frac{p_1}{p_1^2-1} &0 &0 &\dots &\dots &0 \\
        \frac{p_1}{p_1^2-1} &V_{2} &\frac{p_2}{p_2^2-1} &0 &\dots &\dots &0 \\
        0 &\frac{p_2}{p_2^2-1} &V_{3} &\frac{p_3}{p_3^2-1} &\dots &\dots &0 \\
        0 &0 &\frac{p_3}{p_3^2-1} &V_{4} &\dots &\dots &\vdots \\
        \vdots &\vdots &\vdots &\vdots &\ddots & &\vdots \\
        \vdots &\vdots &\vdots &\vdots & &V_{n-1} &\frac{p_{n-1}}{p_{n-1}^2-1} \\
        0 &0 &0 &\dots &\dots &\frac{p_{n-1}}{p_{n-1}^2-1} &\frac{1}{1-p_{n-1}^2}\\
     \end{bmatrix},
\end{equation}
where $V_k:=\frac{1-p_k^2p_{k-1}^2}{(p_k^2-1)(p_{k-1}^2-1)}=\frac
1{1-p_k^2}+ \frac {p_{k-1}^2}{1-p_{k-1}^2}, \ k=2,\dots,n-1$. Obviously,
${\mathbf 1}_{nm}={\mathbf 1}_n\otimes {\mathbf 1}_m$, and in this way
\begin{align*}
M_{\theta}(n,m)={\mathbf 1}_{nm}^{\top} C^{-1}(n,m,r){\mathbf 1}_{nm}=&\big({\mathbf
  1}_n^{\top} \otimes {\mathbf 1}_m^{\top}\big)
\big(P^{-1}(n,r)\otimes Q^{-1}(m,r)\big)\big({\mathbf
  1}_n \otimes {\mathbf 1}_m\big) \\
=&\big({\mathbf
  1}_n^{\top}P^{-1}(n,r){\mathbf 1}_n\big)\big({\mathbf 1}_m^{\top}
Q^{-1}(m,r){\mathbf 1}_m\big).
\end{align*}
Further, by the same arguments as in \citet{bm} we have
\begin{equation}
     \label{eq:eqA.4}
{\mathbf 1}_n^{\top}P^{-1}(n,r){\mathbf
  1}_n=1+\sum_{i=1}^{n-1}\frac{1-p_i}{1+p_i} \qquad \text{and} \qquad
{\mathbf 1}_m^{\top}Q^{-1}(m,r){\mathbf
  1}_m=1+\sum_{j=1}^{m-1}\frac{1-q_j}{1+q_j},
\end{equation}
implying
\begin{equation*}
M_{\theta}(n,m)\!=\!M_{\theta}^{(1)}(n)M_{\theta}^{(2)}(m), \quad
\text{where} \quad
M_{\theta}^{(1)}(n)\!:=\!1\!+\!\sum_{i=1}^{n-1}\frac{1\!-\!p_i}{1\!+\!p_i}, \ \
M_{\theta}^{(2)}(m)\!:=\!1\!+\!\sum_{j=1}^{m-1}\frac{1\!-\!q_j}{1\!+\!q_j}.
\end{equation*}
Now, consider reformulation
\begin{equation*}
M_{\theta}^{(1)}(n)=1+\sum_{i=1}^{n-1}g\big(\alpha
  d_i\big), \qquad M_{\theta}^{(2)}(n)=1+\sum_{j=1}^{m-1}g\big(\beta
  \delta_j\big)\qquad \text{where} \quad
  g(x):=\frac{1-\exp(-x)}{1+\exp(-x)}.
\end{equation*}
As $g(x)$ is a concave function of $x$, by
\citet[Proposition C1, p. 64]{mo}, $M_{\theta}^{(1)}(n)$ and
$M_{\theta}^{(2)}(m)$ are
Schur-concave functions of their arguments $d_i, \ i=1,2,\ldots ,n-1$,
and $\delta_j, \ j=1,2,\ldots ,m-1$, respectively. In this way
$M_{\theta}(n,m)$  attains its maximum at  $d_1=d_2=\ldots =d_{n-1}$ and
$\delta_1=\delta _2=\ldots =\delta_{m-1}$, which completes the proof.
 \proofend

\subsection{Proof of Theorem \ref{Mrn}}
     \label{subs:subsA.2}

By representation \eqref{eq:eqA.1} and the properties of the Kronecker
product we have
\begin{align*}
M_{\alpha}(n,m)&=\frac 12 \tr
\left\{\Bigg(\Big (P^{-1}(n,r)\otimes Q^{-1}(m,r)\Big)\bigg (\frac{\partial
    P(n,r)}{\partial \alpha}\otimes Q(m,r)\bigg)\Bigg)^2 \right\}\\
&= \frac 12 \tr
\left\{\bigg(P^{-1}(n,r)\frac{\partial
    P(n,r)}{\partial \alpha}\bigg)^2 \otimes {\mathcal I}_m \right\}=
\frac m2 \tr
\left\{\bigg(P^{-1}(n,r)\frac{\partial
    P(n,r)}{\partial \alpha}\bigg)^2\right\},
\end{align*}
where ${\mathcal I}_m$ denotes the $m\times m$ unit matrix. Now, the same
ideas that lead to the  proof of
\citet[Theorem 2]{bm} \citep[see also][Proposition 6.1]{baz} imply
the first equation of \eqref{eq:eq4.2}. The form of $M_{\beta}(n,m)$
follows by symmetry. Finally,
\begin{align*}
M_{\alpha,\beta}(n,m)=&\frac 12 \tr
\Bigg\{\big (P^{-1}(n,r)\otimes Q^{-1}(m,r)\big)\bigg (\frac{\partial
    P(n,r)}{\partial \alpha}\otimes Q(m,r)\bigg) \\
   &\phantom{\frac 12 \tr \Bigg\{}\times \big
    (P^{-1}(n,r)\otimes Q^{-1}(m,r)\big) \bigg (P(n,r)\otimes \frac{\partial
    Q(m,r)}{\partial \beta}\bigg) \Bigg\}\\
=& \frac 12 \tr
\left\{\Bigg(\bigg(P^{-1}(n,r)\frac{\partial
    P(n,r)}{\partial \alpha}\bigg) \otimes {\mathcal I}_m\Bigg)
  \Bigg({\mathcal I}_n\otimes \bigg(Q^{-1}(m,r)\frac{\partial
    Q(m,r)}{\partial \beta}\bigg)\Bigg) \right\}\\
=&\frac 12 \tr
\left\{P^{-1}(n,r)\frac{\partial
    P(n,r)}{\partial \alpha}\right\} \tr
\left\{Q^{-1}(m,r)\frac{\partial
    Q(m,r)}{\partial \beta}\right\},
\end{align*}
so the last statement of Theorem \ref{Mrn} follows from \citet[Theorem
3.1]{zba} (see also \citet[Proposition 6.1]{baz}). \proofend

\subsection{Proof of Theorem \ref{covpars}}
       \label{subs:subsA.3}

Consider first the case when we are interested in estimation of
one of the parameters $\alpha$ and $\beta$ and other parameters are
considered as nuisance. According to Remark \ref{rem1}, in this
situation the statement of the theorem directly follows from the
corresponding result for OU processes, see \citet[Theorem 4.2]{zba}

Now, consider the case when both  $\alpha$ and $\beta$ are
unknown. According to \eqref{eq:eq4.1} and  \eqref{eq:eq4.2} the
corresponding objective function to be maximized is
\begin{align}
\label{eq:eqA.5}
\Phi(d_1,&\,\ldots,d_{n-1},\delta_1,\ldots,\delta_{m-1})=\det\big(M_r(n,m)\big)
\\ &=nm
\Bigg (\sum_{i=1}^{n-1}\frac{d_i^2p_i^2(1+p_i^2)}{(1-p_i^2)^2}\Bigg )
\Bigg
(\sum_{j=1}^{m-1}\frac{\delta_j^2q_j^2(1+q_j^2)}{(1-q_j^2)^2}\Bigg )-
4\Bigg (\sum_{i=1}^{n-1}\frac{d_ip_i^2}{1-p_i^2}\Bigg )^2
\Bigg
(\sum_{j=1}^{m-1}\frac{\delta_jq_j^2}{1-q_j^2}\Bigg )^2, \nonumber
\end{align}
which is non-negative, due to Cauchy-Schwartz inequality. Short
calculation shows
\begin{align}
\label{eq:eqA.6}
\Phi(d_1,\ldots,d_{n-1},\delta_1,\ldots,\delta_{m-1})=&\,
\bigg(n\sum_{i=1}^{n-1}g(d_i,
\alpha) \bigg)\left(m\sum_{j=1}^{m-1}g(\delta_j,\beta)-2\bigg(\sum_{j=1}^{m-1}
h(\delta_j,\beta)\bigg)^2\right) \\
&+2\bigg(\sum_{j=1}^{m-1}h(\delta_j,\beta) \bigg)^2
\left(n\sum_{i=1}^{n-1}g(d_i,\alpha)-2\bigg(\sum_{i=1}^{n-1}
h(d_i,\alpha)\bigg)^2\right), \nonumber
\end{align}
where
\begin{equation}
   \label{eq:eqA.7}
g(x,\gamma):=\frac {x^2\big(\exp(2\gamma x)+1\big)}
{\big(\exp(2\gamma x)-1\big)^2} \qquad \text{and} \qquad
h(x,\gamma):=\frac {x}{\exp(2\gamma x)-1}.
\end{equation}
In this way one can consider the two coordinate directions separately.

Since for a given parameter value $\gamma$ both $g(x,\gamma)$
\citep[Theorem 4.2]{zba} and  $h(x,\gamma)$ \citep[Theorem 4.2]{baz}
are convex functions of $x$, according to \citet[Proposition C1, p. 64]{mo}
$$\sum_{i=1}^{n-1}g(d_i, \alpha) \qquad \text{and} \qquad
\bigg(\sum_{j=1}^{m-1}h(\delta_j,\beta) \bigg)^2$$
are Schur-convex functions on $[0,1]^{n-1}$ and  $[0,1]^{m-1}$,
respectively. In this way, they can attain their maxima on the
frontiers of their domains of definition.

Finally, consider the constrained optimum of, e.g.,
\begin{equation*}
\Psi (d_1,\ldots ,d_{n-1}):=n\sum_{i=1}^{n-1}g(d_i,\alpha)-2\bigg(\sum_{i=1}^{n-1}
h(d_i,\alpha)\bigg)^2, \qquad \text{given} \qquad \sum_{i=1}^{n-1}d_i=1.
\end{equation*}
Equating the partial derivatives of the Lagrange function
\begin{equation*}
\Lambda (d_1,\ldots d_{n-1};\lambda):= \Psi (d_1,\ldots
,d_{n-1})+\lambda(d_1+\ldots +d_{n-1}-1)
\end{equation*}
to zero results in equations
\begin{equation*}
n g'(d_k,\alpha)-4\bigg(\sum_{i=1}^{n-1}h(d_i,\alpha)
\bigg)h'(d_k,\alpha)+\lambda =0, \qquad k=1,2,\ldots ,n-1.
\end{equation*}
This means that the  optimum point of $\Psi$ in $[0,1]^{n-1}$ corresponds
to the equidistant design $d_1=d_2=\ldots $ $=d_{n-1}=1/(n-1)$. \proofend

\subsection{Proof of Theorem \ref{freedes1}}
  \label{subs:subsA.4}
Observe first that instead of $\det\big(M_r(n,m)\big)$ given by
\eqref{eq:eq4.3} it suffices to investigate the behaviour of the
function
\begin{equation*}
G(x,y):=\frac{x^2y^2}{({\mathrm
    e}^x-1)^2({\mathrm e}^y-1)^2}\Big(nm ({\mathrm
    e}^x+1)({\mathrm e}^y+1)-4(n-1)(m-1)\Big), \qquad
  x,y\geq 0.
\end{equation*}
Obviously,
\begin{equation*}
\frac{\partial G(x,y)}{\partial x}=\frac{xy^2}{({\mathrm
    e}^x-1)^3({\mathrm e}^y-1)^2}\Big(nm ({\mathrm
  e}^y+1)\big((2-x){\mathrm e}^{2x}-3x{\mathrm
  e}^x-2\big)+8(n-1)(m-1)\big (1-(1-x){\mathrm e}^x\big)\Big),
\end{equation*}
which equals $0$ for non-zero values of $x$ and $y$ if and only if
\begin{equation}
  \label{eq:eqA.8}
\frac{1-(1-x){\mathrm e}^x}{(x-2){\mathrm e}^{2x}+3x{\mathrm
  e}^x+2}=\frac{nm({\mathrm e}^y+1)}{8(n-1)(m-1)}.
\end{equation}
Now, the left-hand side of \eqref{eq:eqA.8} is strictly monotone
decreasing and has a range of $[0,1/2]$. If $nm\geq 2(n-1)(m-1)$ then for
$y>0$ the right-hand side of \eqref{eq:eqA.8} is greater than $1/2$, so in this
case $\frac{\partial G(x,y)}{\partial x}<0$. Finally, if $nm<
2(n-1)(m-1)$ and $y$ is fixed and small enough then the right-hand side of
\eqref{eq:eqA.8} is less than $1/2$, so $\frac{\partial
  G(x,y)}{\partial x}=0$ in a
single point $x$, where $G(x,y)$ takes its maximum. \proofend

\subsection{Calculations for Example \ref{ex:ex2.2}}
  \label{subs:subsA.5}
Decomposition \eqref{eq:eqA.6} of $\det \big(M_r(3,3)\big)$ implies
\begin{align}
\det\big(M(3,3)\big)=&\,\Big[3\big(1+\phi(d,\alpha)+\phi(1-d,\alpha)\big)
\big(g(d,\alpha)+g(1-d,\alpha)\big)\Big] \nonumber \\
&\times
\bigg[\big(1+\phi(\delta,\beta)+\phi(1-\delta,\beta)\big)
\Big(3\big(g(\delta,\beta)
+g(1-\delta,\beta)\big)-2\big(h(\delta,\beta)+h(1-\delta,\beta)\big)^2\Big)
\bigg] \label{eq:eqA.9}\\
&+\bigg[\big(1+\phi(d,\alpha)+\phi(1-d,\alpha)\big)\Big(3\big(g(d,\alpha)
+g(1-d,\alpha)\big)-2\big(h(d,\alpha)+h(1-d,\alpha)\big)^2\Big)
\bigg] \nonumber\\
&\times
\Big[2\big(1+\phi(\delta,\beta)+\phi(1-\delta,\beta)\big)\big(h(\delta,\beta)
+h(1-\delta,\beta)\big)^2\Big], \nonumber
\end{align}
where $g(x,\gamma)$ and $h(x,\gamma)$ are defined by \eqref{eq:eqA.7} and
\begin{equation*}
  \phi(x,\gamma):=\frac{\exp(\gamma x)-1}{\exp(\gamma x)+1}.
\end{equation*}
In this way one can separate $d$ and $\delta$ and  it
suffices to investigate the behaviour of functions
\begin{align*}
\Phi_1(x,\gamma):=&\,\Psi_1(x,\gamma)\Psi_2(x,\gamma), \qquad
\Phi_2(x,\gamma):=\Psi_1(x,\gamma)\big(\Psi_2(x,\gamma)\big)^2,\\
&\Phi_3(x,\gamma):=\Psi_1(x,\gamma)\Big(3\Psi_2(x,\gamma)
-2\big(\Psi_3(x,\gamma)\big)^2\Big),  
\end{align*}
where $x\in[0,1], \ \gamma>0$ and
\begin{align*}
\Psi_1(x,\gamma)&:=1+\phi(x,\gamma)+\phi(1\!-\!x,\gamma), \quad
\Psi_2(x,\gamma):=g(x,\gamma)+g(1\!-\!x,\gamma),\quad
\Psi_3(x,\gamma):=h(x,\gamma)+h(1\!-\!x, \gamma).
\end{align*}
$\Psi_1(x,\gamma),\ \Psi_2(x,\gamma)$ and
$\Psi_3(x,\gamma)$ are symmetric in $x$ on $1/2$ and obviously, the
same property holds for $\Phi_1(x,\gamma),\ \Phi_2(x,\gamma)$ and
$\Phi_3(x,\gamma)$. Further, as
$\frac{\partial \phi (x,\gamma)}{\partial x}$ is strictly monotone
decreasing, while $\frac{\partial g (x,\gamma)}{\partial x}$ and
$\frac{\partial h (x,\gamma)}{\partial x}$ are strictly 
monotone increasing,  $\Psi_1$ is strictly concave, while  $\Psi_2$ and
$\Psi_3$ are strictly convex functions of $x$. 

Consider first $\Phi_1(x,\gamma)$. As
\begin{equation*}
\Psi_1(0,\gamma)\leq
\Psi_1(x,\gamma)\leq \Psi_1(1/2,\gamma) \quad \text{and} \quad
\Psi_2(1/2,\gamma)\leq \Psi_2(x,\gamma)\leq \Psi_2(0,\gamma ), \qquad
x\in [0,1],
\end{equation*}
we have
\begin{equation}
  \label{eq:eqA.10}
\frac{\partial \Phi_1(x,\gamma)}{\partial x}\begin{cases}
\leq{\mathcal Y}(x,\gamma), & \text{if $0<x<1/2$;}\\
\geq{\mathcal Y}(x,\gamma), & \text{if $1/2\leq x<1$,}
\end{cases}
\end{equation}
where
\begin{equation}
  \label{eq:eqA.11}
{\mathcal Y}(x,\gamma):= 
\frac{\partial \Psi_1(x,\gamma)}{\partial
  x}\Psi_2(0,\gamma)+\frac{\partial \Psi_2(x,\gamma)}{\partial
  x}\Psi_1(0,\gamma)=\Upsilon(x,\gamma)-\Upsilon(1-x,\gamma),
\end{equation}
with
\begin{equation*}
\Upsilon(x,\gamma): =\frac{\partial \phi(x,\gamma)}{\partial
  x}\frac{{\mathrm e}^{2\gamma}\!+\!1}{({\mathrm
    e}^{2\gamma}\!-\!1)^2}+\frac{\partial g(x,\gamma)}{\partial
  x}\frac{2{\mathrm e}^{\gamma}}{{\mathrm e}^{\gamma}\!+\!1}=
\frac{2\gamma {\mathrm e}^{\gamma x}({\mathrm e}^{2\gamma}\!+\!1)}{({\mathrm
    e}^{\gamma x}\!+\!1)^2({\mathrm
    e}^{2\gamma}\!-\!1)^2}-\frac{4x{\mathrm e}^{\gamma}(3\gamma x
  {\mathrm e}^{2\gamma x}\!-\!{\mathrm e}^{4\gamma x}\!+\!\gamma x {\mathrm
    e}^{4\gamma x}\!+\!1)}{({\mathrm e}^{2\gamma x}\!-\!1)^3 ({\mathrm
    e}^{\gamma}\!+\!1)}. 
\end{equation*}
Further, let
\begin{equation*}
\frac{\partial \Upsilon (x,\gamma)}{\partial
  x}=\frac {\Upsilon ^{(1)} (x,\gamma)}{\Upsilon ^{(2)} (x,\gamma)},
\end{equation*}
where for $x>0$ the denominator $\Upsilon ^{(2)}
(x,\gamma)=({\mathrm e}^{2\gamma x}-1)^4({\mathrm e}^{2\gamma}-1)^2$ is  
obviously positive, while the numerator can be written as
\begin{align*}
\Upsilon ^{(1)} (x,\gamma)=&\,4{\mathrm e}^{\gamma}({\mathrm
  e}^{2\gamma}\!-\!1)({\mathrm e}^{\gamma}\!-\!1)\Big({\mathrm
  e}^{6\gamma x}\big(2\gamma^2x^2\!-\!4\gamma x\!+\!1\big)\!+\!{\mathrm
  e}^{4\gamma x}\big(16\gamma^2x^2\!-\!8\gamma x\!-\!1\big)\!+\!{\mathrm
  e}^{2\gamma x}\big(6\gamma^2x^2\!+\!12\gamma x\!-\!1\big)\!+\!1\Big)
\\ &-2\gamma^2({\mathrm e}^{2\gamma}\!+\!1){\mathrm e}^{\gamma
  x}({\mathrm e}^{\gamma x}\!-\!1)^5({\mathrm e}^{\gamma x}\!+\!1).
\end{align*}
If $x\in[0,1]$ then by inequality
\begin{equation*}
2{\mathrm e}^{\gamma}({\mathrm
  e}^{2\gamma}\!-\!1)({\mathrm e}^{\gamma}\!-\!1)>\gamma^2({\mathrm
  e}^{2\gamma}\!+\!1), \qquad \gamma>0, 
\end{equation*}
we have
\begin{equation}
   \label{eq:eqA.12}
\Upsilon ^{(1)} (x,\gamma)\geq \gamma^2{\mathrm e}^{\gamma}({\mathrm
  e}^{2\gamma}\!+\!1)S(\gamma x),
\end{equation}
where
\begin{equation*}
S(y):={\mathrm e}^{6y}(y^2-2y)+2{\mathrm e}^{5y}+{\mathrm
  e}^{4y}(8y^2-4y-3)+{\mathrm e}^{2y}(3y^2+6y+2)-2{\mathrm e}^{y}+1. 
\end{equation*}
Short calculation shows that $S(y)$ is positive if $y>0$, which
together with \eqref{eq:eqA.12} implies the positivity of $\Upsilon
^{(1)} (x,\gamma)$ for $0<x\leq 1$. Thus, $\Upsilon (x,\gamma)$ is
strictly monotone increasing, so using \eqref{eq:eqA.11} one can easily
see that ${\mathcal Y}(x,\gamma)<0$ if $x<1/2$.
Now, \eqref{eq:eqA.10} implies that $\Phi_1(x,\gamma)$ has 
a single global minimum at $1/2$, while its maximum is reached at $0$
and $1$. In a similar way one can 
verify that $\Phi_2(x,\gamma)$ and $\Phi_3(x,\gamma)$ have the same
behaviour, and since all coefficients in \eqref{eq:eqA.9} are
non-negative, this completes the proof. \proofend

\subsection{Proof of Theorem \ref{freedes2}}
   \label{subs:subsA.6}
Similarly to the proof of Theorem \ref{freedes1}, instead of
$\det\big(M(n,n)\big)$ given by \eqref{eq:eq4.5} one can consider
function
\begin{equation*}
G(x,y):=\frac{x^2y^2\big(n({\mathrm e}^{x}-1)+2\big)\big(n({\mathrm
    e}^{y}-1)+2\big) }{\big({\mathrm e}^{2x}-1\big)^2\big({\mathrm
    e}^{2y}-1\big)^2\big({\mathrm e}^{x}+1\big)\big({\mathrm e}^{y}+1\big)}
\Big(n^2 \big({\mathrm
    e}^{2x}+1\big)\big({\mathrm e}^{2y}+1\big)-4(n-1)^2\Big),
  \qquad x,y\geq 0.
\end{equation*}
Short calculation shows
\begin{equation*}
\frac{\partial G(x,y)}{\partial x}=\frac{2xy^2\big(n({\mathrm
    e}^{y}-1)+2\big)}{\big({\mathrm
    e}^{2x}-1\big)^3\big({\mathrm
    e}^{2y}-1\big)^2\big({\mathrm e}^{x}+1\big)\big({\mathrm
    e}^{y}+1\big)} \Big(n^2 \big({\mathrm
  e}^{2y}+1\big)g_1(x,n)-4(n-1)^2g_2(x,n)\Big),
\end{equation*}
where
$g_1(x,n)$ and $g_2(x,n)$ are defined by \eqref{eq:eq4.7}.
Hence, the extremal points of $G(x,y)$ should solve
\begin{equation*}
n^2\big({\mathrm e}^{2y}+1\big)g_1(x,n)=4(n-1)^2g_2(x,n), \qquad
n^2\big({\mathrm e}^{2x}+1\big)g_1(y,n)=4(n-1)^2g_2(y,n),
\end{equation*}
which proves \eqref{eq:eq4.6}.

Assume first $n=2$. In this case $g_2(x,n)/g_1(x,n)$ is strictly
monotone decreasing and has a range of $[0,3/2]$, while
$n^2\big({\mathrm e}^{2y}+1\big)/\big(4(n-1)^2\big)>3/2$, implying
$\frac{\partial G(x,y)}{\partial x}<0$.

Now, let us fix $y>0$ and assume $n\geq 3$. In this case
\begin{equation*}
\lim_{x\searrow 0}\frac{\partial G(x,y)}{\partial
  x}=\frac{(n\!-\!1)^2y^2\big(n({\mathrm e}^{y}\!-\!1)\!+\!2\big)}{4\big({\mathrm
    e}^{2y}-1\big)^2\big({\mathrm e}^{y}+1\big)}\big (n^2 ({\mathrm
  e}^{2y}-1)(n-3)+ 2(4n^2-11n+5)\big)>0 \quad \text{and} \quad
\lim_{x\to\infty}\frac{\partial G(x,y)}{\partial x}=0,
\end{equation*}
so $G(x,y)\geq 0$ should have a global maximum  at some $x>0$. The
same result can be proved if we fix $x>0$ and
consider $G(x,y)$ as a function of $y$. This means that if $n\geq 3$
then $G(x,y)$ reaches its global maximum at a point with non-zero
coordinates, which completes the proof.

\subsection{Proof of Theorem \ref{IMSPE}}
        \label{subs:subsA.7}

Observe first, that the product structure of elements of
$R(x_1,x_2)$
implies that $R(x_1,x_2)=R_1(x_1)\otimes
R_2(x_2)$ with $R_1(x_1)=(\varrho_{1,1},
\varrho_{1,2}, \ldots ,\varrho_{1,n})^{\top}$   and
$R_2(x_2)=(\varrho_{2,1}, \varrho_{2,2}, \ldots
,\varrho_{2,m})^{\top}$, where to shorten our formulae instead of
$\varrho_1(x_1,s_i)$ and $\varrho_2(x_2,t_j)$ we use simply
$\varrho_{1,i}$ and $\varrho_{2,j}$, respectively, $i=1,2,\ldots ,
n, \ j=1,2, \ldots ,m$.

Consider first $\mspe\big (\widehat
Y(x_1,x_2)\big)$ given by \eqref{eq:eq6.1}. Using matrix algebraic
calculations \citep[see, e.g.,][]{bss}, decomposition of $R(x_1,x_2)$
and \eqref{eq:eqA.2}, one can easily show
\begin{align}
\mspe\big (&\,\widehat Y(x_1,x_2)\big) \nonumber \\
=&\,\sigma^2\Big[1-R^{\top}(x_1,x_2)C^{-1}(n,m,r)
R(x_1,x_2)+M_{\theta}^{-1}(n,m)\big(1-R^{\top}(x_1,x_2)C^{-1}(n,m,r)
{\mathbf 1}_{nm}\big)^2\Big]  \label{eq:eqA.13} \\
=&\,\sigma^2\bigg[1-\big(R_1^{\top}(x_1)P^{-1}(n,r)R_1(x_1)\big)
\big(R_2^{\top}(x_2)Q^{-1}(m,r) R_2(x_2)\big)\! \nonumber\\
&\phantom{==}+M_{\theta}^{-1}(n,m)\Big(1-\big(R_1^{\top}(x_1)P^{-1}(n,r)
{\mathbf 1}_n\big)\big(R_2^{\top}(x_2)Q^{-1}(m,r) {\mathbf
  1}_m\big)\Big)^2\bigg],
 \nonumber
\end{align}
which implies \eqref{eq:eq6.2}.

Further, according to the definition of IMSPE criterion, we can write
\begin{equation*}
\imspe\big (\widehat Y\big)=1-A^{(1)}_nA^{(2)}_m+\Big(1+\sum_{i=1}^{n-1}
\frac{1-p_i}{1+p_i}\Big)^{-1}\Big(1+\sum_{j=1}^{m-1}
\frac{1-q_j}{1+q_j}\Big)^{-1}\Big(1-2B^{(1)}_nB^{(2)}_m+D^{(1)}_nD^{(2)}_m\Big),
\end{equation*}
where
\begin{alignat*}{3}
A^{(1)}_n:=&\,
\tr \left[P^{-1}(n,r){\mathcal R}_1\right], \quad &B^{(1)}_n:={\mathbf
  1}_n^{\top}P^{-1}(n,r){\mathcal W}_1, \quad  &D^{(1)}_n:={\mathbf
  1}_n^{\top}P^{-1}(n,r){\mathcal R}_1P^{-1}(n,r){\mathbf 1}_n, \\
A^{(2)}_m:=&\,\tr
\left[Q^{-1}(m,r){\mathcal R}_2\right], \quad &B^{(2)}_m:={\mathbf
  1}_m^{\top}Q^{-1}(m,r){\mathcal W}_2 \quad  &D^{(2)}_m:={\mathbf
  1}_m^{\top}Q^{-1}(m,r){\mathcal R}_2Q^{-1}(m,r){\mathbf 1}_m ,
\end{alignat*}
with
\begin{equation*}
{\mathcal W}_s=\big\{\omega_{s,i}\big\}:=\int\limits_0^1R_s(x)\,
{\mathrm d}x \qquad \text {and} \qquad {\mathcal
  R}_s=\big\{
R_{s,i,j}\big\}:=\int\limits_0^1R_s(x)R_s^{\top}(x)\,{\mathrm d}x,
\qquad s=1,2.
\end{equation*}
Obviously,
\begin{align}
\omega_{1,i}=&\,\frac 1{\alpha}\Big[2-{\mathrm e}^{-\alpha s_i}-
{\mathrm e}^{-\alpha (1-s_i)}\Big], \quad \quad
\omega_{2,i}=\frac 1{\beta}\Big[2-{\mathrm e}^{-\beta t_i}-
{\mathrm e}^{-\beta (1-t_i)}\Big],  \nonumber\\
R_{1,i,j}=&\,\frac 1{2\alpha}\Big(2{\mathrm e}^{-\alpha
  |s_i-s_j|}-{\mathrm e}^{-\alpha (s_i+s_j)}-{\mathrm e}^{-\alpha
  (2-s_i-s_j)}\Big)+|s_i-s_j|{\mathrm e}^{-\alpha |s_i-s_j|}, \label{eq:eqA.14}\\
R_{2,i,j}=&\,\frac 1{2\beta}\Big(2{\mathrm e}^{-\beta
  |t_i-t_j|}-{\mathrm e}^{-\beta (t_i+t_j)}-{\mathrm e}^{-\beta
  (2-t_i-t_j)}\Big)+|t_i-t_j|{\mathrm e}^{-\beta |t_i-t_j|}. \nonumber
\end{align}
Now, extracting, e.g., the expressions for $A^{(1)}_n, \ B^{(1)}_n$ and
$D^{(1)}_n$ we obtain
\begin{align*}
A^{(1)}_n=&\,R_{1,n,n}+\sum_{i=1}^{n-1}\frac
{R_{1,i,i}-2R_{1,i+1,i}p_i+R_{1,i+1,i+1}p_i^2}{1-p_i^2}, \qquad \qquad
B^{(1)}_n=\omega_{1,n}+\sum_{i=1}^{n-1}\frac
{\omega_{1,i}-\omega_{1,i+1}p_i}{1+p_i}, \\
D^{(1)}_n=&\,R_{1,n,n}+2\sum_{i=1}^{n-1}\frac{R_{1,n,i}-
  R_{1,n,i+1}p_i}{1+p_i}
+\sum_{i=1}^{n-1}\sum_{j=1}^{n-1}
\frac{R_{1,i,j}-R_{1,i+1,j}p_i -R_{1,i,j+1}p_j+R_{1,i+1,j+1}p_ip_j}{(1+p_i)(1+p_j)},
\end{align*}
and long but straightforward calculations using \eqref{eq:eqA.14}
yield
\begin{equation*}
A^{(1)}_n=\frac {n-1}{\alpha}-2\sum_{i=1}^{n-1}\frac
{d_ip_i^2}{1-p_i^2}, \qquad B^{(1)}_n=\frac 2{\alpha}\sum_{i=1}^{n-1}\frac
{1-p_i}{1+p_i}, \qquad
D^{(1)}_n=\sum_{i=1}^{n-1}\frac
{1-p_i^2+2\alpha d_ip_i}{\alpha (1+p_i)^2}.
\end{equation*}
The closed forms of $A^{(2)}_m, \ B^{(2)}_m$ and $D^{(2)}_m$ can be
derived in the same way.

Obviously, $\imspe\big (\widehat Y\big)$ is permutation invariant with
respect to both $d_1,d_2, \ldots ,d_{n-1}$ and $\delta_1, \delta_2,
\ldots, \delta_{m-1}$. Now, fix, e.g., $\delta_1, \delta_2, \ldots,
\delta_{m-1}$ and consider the partial derivatives
\begin{align}
 \frac{\partial \imspe\big (\widehat Y\big)}{\partial d_i}=&
2\frac{\partial h(d_i,\alpha)}{\partial d}\bigg(\frac
{m-1}{\beta}-2H_m(\boldsymbol\delta,\beta)\bigg)+\frac
4{\alpha} \bigg (\frac{\partial \varphi (d_i,\alpha)}{\partial
  d}-1\bigg)\frac{\Psi_m(\boldsymbol\delta,\beta)-2\big(\Phi_m(
  \boldsymbol\delta,\beta)-1\big)/\beta}
{\Phi_n(\boldsymbol d,\alpha)\Phi_m(\boldsymbol\delta,\beta)}  \nonumber
\\
& +\bigg(\frac{\partial \psi (d_i,\alpha)}{\partial d}-\frac
4{\alpha} \Big (\frac{\partial \varphi (d_i,\alpha)}{\partial
  d}-1\Big)\bigg)
\frac{\Psi_m(\boldsymbol\delta,\beta)}
{\Phi_n(\boldsymbol d,\alpha)\Phi_m(\boldsymbol\delta,\beta)}
 \label{eq:eqA.15} \\
&- \frac{\partial \varphi (d_i,\alpha)}{\partial d}\frac{1-
8\big(\Phi_m(\boldsymbol d,\alpha)-1\big)
\big(\Phi_m(\boldsymbol\delta,\beta)-1\big)/(\alpha\beta)
+\Psi_n(\boldsymbol d,\alpha)\Psi_m(\boldsymbol\delta,\beta)}{
\Phi_n(\boldsymbol d,\alpha)^2\Phi_m(\boldsymbol\delta,\beta)}, \nonumber
\end{align}
where
\begin{equation*}
h(x,\gamma):=\frac {x}{\exp(2\gamma x)-1}, \qquad
\varphi (x,\gamma):=x+\frac {\exp(\gamma x)-1}{\exp(\gamma x)+1},
\qquad
\psi (x,\gamma):=\frac {\exp(2\gamma x)-1+2\gamma x \exp(\gamma x)
}{\gamma (\exp(\gamma x)+1)^2},
\end{equation*}
and for $x_1,x_2,\ldots, x_{n-1}$ define
\begin{equation*}
H_n(\boldsymbol x,\gamma):=\sum_{i=1}^{n-1}h(x_i,\gamma), \qquad
\Phi_n(\boldsymbol x,\gamma):=\sum_{i=1}^{n-1}\varphi(x_i,\gamma), \qquad
\Psi_n(\boldsymbol x,\gamma):=\sum_{i=1}^{n-1}\psi(x_i,\gamma).
\end{equation*}

Short calculation shows \citep[see, e.g.,][]{baz} that on the $[0,1]$ interval
$\varphi(x,\gamma)$ is concave, while $h(x,\gamma)$ and
$\psi(x,\gamma)-4\big(\varphi(x,\gamma)-x\big)/\gamma $ are convex
functions of $x$.  Further,
for $x_i\geq 0, \ i=1,2, \ldots ,n-1$, we have $\Psi_n(\boldsymbol
x,\gamma)\geq 0$, inequality
$\exp(x)-1\geq x, \ x\in {\mathbb R}$, implies $2H_n(\boldsymbol x,\gamma)\leq
(n-1)/\gamma$, and if in addition we assume  $\sum_{i=1}^{n-1}x_i=1$, then
$\Phi_n(\boldsymbol x,\gamma)\geq 1$ and $\gamma\Psi_n(\boldsymbol
x,\gamma)\leq 2\Phi_n(\boldsymbol x,\gamma)-2$ also hold. Finally,
representation
\eqref{eq:eqA.13} of the $\mspe$ implies that the
numerator of the fraction in the last term \eqref{eq:eqA.15} is also
non-negative, so
$\frac{\partial \imspe\big (\widehat Y\big)}{\partial d_i}$ is
monotone increasing in $d_i$. Hence, for all fixed $\delta_1, \delta_2, \ldots,
\delta_{m-1}$ function $\imspe\big (\widehat Y\big)$ is
Schur convex \citep[see, e.g.,][Theorem A.4, p. 57]{mo}, so it attains its
minimum at $d_i=1/(n-1), \ i=1,2,\ldots ,n-1$. An analogous result can
be derived if we fix $d_1,d_2, \ldots ,d_n$ and consider
$\imspe\big (\widehat Y\big)$ as a function of  $\delta_1, \delta_2, \ldots,
\delta_{m-1}$, which together with the previous statement implies the
optimality of the directionally equidistant design.
\proofend

\subsection{Proof of Theorem \ref{entropy}}
  \label{subs:subsA.8}

Using decomposition \eqref{eq:eqA.1} and the properties of the
Kronecker product one has
\begin{equation*}
\det C(n,m,r)= \big(\det P(n,r)\big)^m \big(\det Q(m,r)\big)^n,
\end{equation*}
hence
\begin{equation*}
\ent ({\mathbf Y})=\frac {nm}2\big(1+\ln (2\pi \sigma^2)\big)+\frac
m2\ln\det P(n,r)+\frac n2\ln\det Q(m,r).
\end{equation*}
The special forms of matrices $P(n,r)$ and $Q(m,r)$ imply \citep[see,
e.g.,][Lemma 3.1]{baz} that
\begin{equation*}
\det P(n,r) =\prod_{i=1}^{n-1} (1-p_i^2) \qquad \text{and} \qquad \det
Q(m,r) =\prod_{j=1}^{m-1} (1-q_j^2),
\end{equation*}
which proves \eqref{eq:eq6.1}.

In order to find the optimal design one has to
find the constrained maximum of
\begin{equation*}
F\big(p_1,\ldots ,p_{n-1},q_1,\ldots ,q_{m-1}\big):=\frac m2 \sum
_{i=1}^{n-1}\ln \big(1-p_i^2\big)+\frac n2\sum _{j=1}^{m-1}\ln
\big(1-q_j^2\big)
\end{equation*}
under conditions
\begin{equation*}
\sum_{i=1}^{n-1} \ln p_i= -\alpha \qquad \text{and} \qquad
\sum_{j=1}^{m-1} \ln q_i= -\beta.
\end{equation*}
By analyzing the first partial derivatives and the Hessian of the
Lagrange function
\begin{align*}
\Lambda\big(p_1,\ldots ,p_{n-1},q_1,\ldots ,q_{m-1};\lambda,\mu
\big):=&\,\frac m2 \sum
_{i=1}^{n-1}\ln \big(1-p_i^2\big)+\frac n2\sum _{j=1}^{m-1}\ln
\big(1-q_j^2\big)\\
&+\lambda \Bigg(
\sum_{i=1}^{n-1} \ln p_i +\alpha\Bigg)+\mu\Bigg(
\sum_{j=1}^{m-1} \ln q_j +\beta\Bigg)
\end{align*}
one can easily see that the maximum is reached when $p_1=p_2=\ldots =p_{n-1}$
and  $q_1=q_2=\ldots =q_{m-1}$,  which completes the proof. \proofend

\subsection{Proof of Theorem \ref{arhfish}}
  \label{subs:subsA.9}
Since
\begin{equation*}
\frac{\partial \eta(s,t;\mu,B)}{\partial B}=-\frac
1{t^{\mu+1}}{\mathrm e}^{-B/t},
\end{equation*}
vector $F(n,m,B)$ can be decomposed as $F(n,m,B)={\mathbf 1}_n\otimes
K(m,B)$, where $K(m,B):=(\kappa_1,\kappa_2,\dots ,\kappa
_m)^{\top}$. Hence, decomposition \eqref{eq:eqA.2} and the properties of
the Kronecker product imply
\begin{equation*}
M_B(n,m)=\big({\mathbf 1}_n^{\top}P^{-1}(n,r){\mathbf 1}_n\big)
\big(K^{\top}(m,B)Q^{-1}(m,r)K(m,B)\big).
\end{equation*}
Using the same calculations as in the proof of \eqref{eq:eq6.2} one
can derive
\begin{equation*}
K^{\top}(m,B)Q^{-1}(m,r)K(m,B)=\kappa_m^2
+\sum_{j=1}^{m-1}\frac{(\kappa_j-\kappa_{j+1}q_j)^2}{1-q_j^2},
\end{equation*}
which together with \eqref{eq:eqA.4} implies
\eqref{arhfishform}. \proofend

\subsection{Proof of Theorem \ref{arhfishmu}}
  \label{subs:subsA.10}
Having a look at the partial derivatives of $\eta(s,t;\mu,B)$ with
respect to $\mu$ and $B$ one can easily see that
$G(n,m,\mu,B)={\mathbf 1}_n\otimes \Lambda(m,\mu,B)$, where
\begin{equation*}
\Lambda(m,\mu,B):=\begin{bmatrix}
\lambda_1&\lambda_2&\dots &\lambda _m\\
\kappa_1&\kappa_2&\dots &\kappa _m
\end{bmatrix}^{\top}.
\end{equation*}
Hence, \eqref{arhfishmuform} can be proved in the same way as
\eqref{arhfishform} has been. \proofend

\subsection{Calculations for Example \ref{ex:ex3.2}}
  \label{subs:subsA.11}
Consider first ${\mathcal M}(2,2)$ as a function of $d$. Obviously,
\begin{equation*}
{\mathcal M}(2,2)=
\frac{8\exp(-2B/\delta)\exp(-2\beta\delta)}{\delta^{2\mu}\big(1-\exp(-2\beta
  \delta)\big)^3}\,Q(d,\delta), \qquad d,\delta\geq 0,
\end{equation*}
where
\begin{equation*}
Q(d,\delta):=
\frac{d^2\exp(-2\alpha
  d)\big(1+\exp(-2\alpha d)+q^2\big)}{\big( 1-\exp(-2\alpha d)\big)^2
  \big(1+\exp(-\alpha d)\big)}, \qquad \text{with} \qquad
q:=\exp(-\beta\delta).
\end{equation*}
Short calculation shows
\begin{equation*}
\frac{\partial Q(d,\delta)}{\partial d}=\frac {-d\exp(\alpha d)}{\big(
  \exp(\alpha d)-1\big)^3 \big(\exp(\alpha d)+1\big)^4}\, S(\alpha d),
\end{equation*}
where
\begin{equation*}
S(x):=x+2-x{\mathrm e}^x+\big(2q^2+7x+3q^2x\big){\mathrm
  e}^{2x}-x\big(1+q^2\big) {\mathrm
  e}^{3x}+2(x-1)(1+q^2){\mathrm e}^{4x}, \qquad x\geq 0.
\end{equation*}
First, let $x\geq 2$ implying $2(x-1)\geq x$, so
\begin{equation*}
S(x)\geq x+2 + x\big({\mathrm e}^{2x}-{\mathrm
  e}^x\big)+x\big(1+q^2\big)\big ({\mathrm
  e}^{4x}- {\mathrm e}^{3x}\big)> 0.
\end{equation*}
Further, for $1<x<2$ we have
\begin{equation*}
S(x)\geq \big(1+q^2\big){\mathrm e}^{2x}\big({\mathrm
  e}^{2x}(x-1)+3x+1-{\mathrm e}^x\big) \geq 0.
\end{equation*}
Finally, consider decomposition $S(x)=\big(1+q^2\big)S_1(x)+S_2(x)$,
where
\begin{equation*}
S_1(x):=2{\mathrm e}^{4x}(x-1)-x{\mathrm
  e}^{3x}+(3x+2){\mathrm e}^{2x} \quad \text{and} \quad
S_2(x):=(4x-2){\mathrm e}^{2x}-x{\mathrm e}^x+x+2.
\end{equation*}
If $0<x\leq 1$ then $S_2(x)\geq 0$ and $2S_1(x)+S_2(x)\geq 0$, which
together with $0<q\leq 1$ imply that on this interval $S(x)$ is
non-negative, too. Hence, for $d\geq 0$ we have $\partial
Q(d,\delta)/\partial d\leq 0$, so $\mathcal M(2,2)$ is decreasing in $d$.

Now, let us investigate decomposition
\begin{equation*}
{\mathcal M}(2,2)=\frac{8d^2\exp(-2\alpha d)}{\big( 1-\exp(-2\alpha d)\big)^2
  \big(1+\exp(-\alpha d)\big)}\, R(d,\delta), \qquad d,\delta\geq 0,
\end{equation*}
where
\begin{equation*}
R(d,\delta)=\frac{\exp(-2B/\delta-2\beta\delta)\big(1+\exp(-2\beta
\delta)+p^2\big)}{\delta^{2\mu}\big(1-\exp(-2\beta \delta)\big)^3}, \qquad
\text{with} \qquad p:=\exp(-\alpha d).
\end{equation*}
Taking the partial derivative of $R$ with respect to $\delta$, after
some calculations we obtain
\begin{equation*}
\frac{\partial R(d,\delta)}{\partial \delta}=
\frac{-2\exp(4\beta\delta)}{\delta
  ^{2\mu+2}\exp(2B/\delta+2\beta\delta) \big(\exp(2\beta
  \delta)-1\big)^4}\,U(\delta),
\end{equation*}
where
\begin{equation*}
U(\delta):=\beta \delta^2-\mu\delta +B+{\mathrm
  e}^{2\beta\delta}\big(2\beta(2+p^2)\delta^2+(B-\mu\delta)p^2 \big)+{\mathrm
  e}^{4\beta\delta}(1+p^2)\big(\beta \delta^2 +\mu\delta -B\big).
\end{equation*}
If $\delta \ne B/\mu$ then equation $U(\delta)=0$ is equivalent to
$V(\delta)=W(\delta)$, where
\begin{equation*}
V(\delta):=\frac{1+p^2-p^2{\mathrm e}^{-2\beta\delta}-{\mathrm
    e}^{-4\beta\delta}}{1+p^2+2(2+p^2){\mathrm e}^{-2\beta\delta}+{\mathrm
    e}^{-4\beta\delta}} \qquad \text{and} \qquad W(\delta):=\frac
{\beta\delta^2}{B -\mu\delta}.
\end{equation*}

Now, let us fix a value $0\leq p<1$. First, consider the function
$V(\delta)$, where
without loss of generality we may assume $\beta=1$. 
One can easily show that $V(\delta)$ is monotone increasing,
$\lim_{\delta\searrow 0}\!V(\delta)=0,\
\lim_{\delta\to\infty}\!V(\delta)=1$ and $\lim_{\delta\searrow
  0}V'(\delta)>0$. Further, $V''(\delta)=-4{\mathrm e}^{2\delta}
\big({\mathcal V}_1(\delta) -{\mathcal
  V}_2(\delta)\big)/\big({\mathcal V}_3(\delta)\big)^3$ with
\begin{align*}
{\mathcal V}_1(\delta)&:=(3p^6+10p^4+11p^2+4){\mathrm
  e}^{8\delta}+(2p^4+4p^2+8){\mathrm e}^{2\delta}>0, \\
{\mathcal V}_2(\delta)&:=(6p^6+18p^4+20p^2+8){\mathrm
  e}^{6\delta}+(6p^4+6p^2){\mathrm e}^{4\delta}+p^2+4>0,\\
{\mathcal V}_3(\delta)&:=(p^2+1){\mathrm
  e}^{4\delta}+(2p^2+4){\mathrm e}^{2\delta}+1>0.
\end{align*}
As both ${\mathcal V}_1(\delta)$ and ${\mathcal V}_2(\delta)$ are
strictly monotone increasing and convex functions,
$\lim_{\delta\searrow 0}\big({\mathcal V}_1(\delta)-{\mathcal
V}_2(\delta)\big)=-3p^6-12p^4-12p^2<0$ and $\lim_{\delta\to \infty}\big({\mathcal
V}_1(\delta)-{\mathcal 
V}_2(\delta)\big)=\infty$, equation ${\mathcal
V}_1(\delta)={\mathcal V}_2(\delta)$ has a single positive root
$\tilde\delta$. This implies that $V(\delta)$ is convex
if $0<\delta<\tilde\delta$ and concave if $\delta>\tilde\delta$.

Concerning the behaviour of $W(\delta)$, assume first $\mu>0$. In
this case $\lim_{\delta\searrow 0}\!W(\delta)=0,\
\lim_{\delta\nearrow B/\mu}\!W(\delta)=\infty$ and $\lim_{\delta\searrow
  B/\mu}W(\delta)=-\infty ,\
\lim_{\delta\to\infty}W(\delta)=-\infty$. Further, for $\delta>B/\mu$ function
$W(\delta)$ has a global maximum at $\delta^*:=2B/\mu$ with
$W(\delta^*)<0$, so on this interval $W(\delta)<V(\delta)$. Finally, if
$0<\delta<B/\mu$ then $W(\delta)$ is strictly monotone increasing and
convex with $\lim_{\delta\searrow 0}W'(\delta)=0$. Hence, for $\mu>0$ equation
$V(\delta)=W(\delta)$ has a single solution which is in the interval
$]0,B/\mu[$. Obviously, if $\mu\leq 0$ then $W(\delta)$ in strictly
monotone increasing and convex on its whole domain of definition. In
this case $\lim_{\delta\searrow 0}W(\delta)=0,\
\lim_{\delta\to\infty}W(\delta)=\infty$ and $\lim_{\delta\searrow
  0}W'(\delta)=0$, so again, the graphs of $V(\delta)$ and $W(\delta)$
intersect in a single point.

As $U(B/\mu)\ne 0$, the above reasoning implies that for any fixed $d$
function $R(d,\delta)$ (and in this way ${\mathcal M}(2,2)$) has a
single extremal point in $\delta$. Since
$\lim_{\delta\searrow  0}R(d,\delta)=0, \ \lim_{\delta\to\infty}R(d,\delta)=0, \
  R(d,\delta)\geq 0$ and $R(d,\delta )\not\equiv 0$, this extremal
  point should be a maximum. \proofend

\end{appendix}

\end{document}